\documentclass[10pt]{article}%
\usepackage{amsfonts}
\usepackage{graphicx}
\usepackage{amsmath}
\usepackage{amssymb}%
\usepackage[a4paper,bottom=3.5cm,top=2cm,left=2.5cm,right=2.5cm]{geometry}
\providecommand{\U}[1]{\protect\rule{.1in}{.1in}}

\newcommand{\R}{\mathbb{R}}

\DeclareMathOperator{\Var}{Var}

\numberwithin{equation}{section}

\newcommand{\Monm}{\Mon^c}
\DeclareMathOperator{\Mon}{Mon}

\renewcommand{\P}{\ensuremath{\mathbb {P}}}
\newcommand{\E}{\mathbb{E}}
\newcommand{\LL}{ { \mathbb L}}
\newcommand{\p}{ { \mathbb P} }

\newcommand{\F}{\mathcal{F}}
\newcommand{\G}{{\mathcal G}}
\newcommand{\h}{{\mathcal H}}

\newcommand{\I}{{\bf 1 }}
\newcommand{\as}{\mbox{$\P$-a.s.}}
\newcommand{\nuas}{\mbox{$\nu$-a.s.}}

\newcommand\beq{\begin{equation}}
\newcommand\eeq{\end{equation}}
\def\E{\mathbb{E}}

\def\R{\mathbb{R}}

\def\F{{\cal F}}

\newtheorem{thm}{Theorem}[section]
\newtheorem{lma}[thm]{Lemma}
\newtheorem{rmk}[thm]{Remark}
\newtheorem{prop}[thm]{Proposition}
\newtheorem{cor}[thm]{Corollary}

\newtheorem{defn}[thm]{Definition}

\newtheorem{claim}{Claim}

\begin{document}

\begin{center}
{\Large \bf Strong invariance principles with rate  for "reverse" martingales 
and applications.}

\bigskip

\bigskip Christophe Cuny$^{a}$ and Florence Merlev\`{e}de$^{b}$
\end{center}

$^{a}$ Laboratoire MAS, Ecole Centrale de Paris, Grande Voie des Vignes, 92295 Chatenay-Malabry cedex, FRANCE. E-mail: christophe.cuny@ecp.fr

$^{b}$ Universit\'{e} Paris Est, LAMA, CNRS UMR 8050, B\^{a}timent Copernic, 5 Boulevard Descartes, 77435 Champs-Sur-Marne,
FRANCE. E-mail: florence.merlevede@univ-mlv.fr

\begin{abstract}

In this paper, we obtain almost sure invariance principles with rate
of order $n^{1/p}\log^\beta n$, $2< p\le 4$,
for sums associated to a sequence of reverse
martingale differences.  Then, we apply those
results to obtain similar conclusions in the context of some
non-invertible dynamical systems. For instance we treat several classes of uniformly expanding maps of the interval (for possibly unbounded functions).
A general result for $\phi$-dependent sequences is obtained in the course.

\end{abstract}

\section{Introduction}
 \setcounter{equation}{0}
Let $T$ be a map from $[0,1]$ to $[0,1]$ preserving a probability $\nu$  on $[0,1]$.  Let $f$ be a measurable function such that $\nu (|f|) < \infty$.  Define  $S_n(f) = \sum_{i=0}^{n-1} (f \circ T^i - \nu (f))$. According to the Birkhoff-Khintchine ergodic theorem, $n^{-1}S_n(f)$ satisfies the strong law of large numbers.
One can go further in the study of the statistical properties of $S_n(f)$ and recently there have been extensive researches to study the almost sure invariance principle (ASIP) for $S_n(f)$. Roughly speaking, such a result ensures that the trajectories of a process can be matched with
the trajectories of a Brownian motion in such a way that almost surely the error
between the trajectories is negligible compared to the size of the trajectory. This result is more or less precise depending on the error term one obtains. In this paper, for some classes of uniformly expanding maps, we give conditions on $f$ ensuring that there exists a sequence
of independent identically distributed (iid) Gaussian random variables $(Z_i)_{i\geq 1}$ such that
\begin{equation}\label{ASIPp}
  \sup_{1 \leq k \leq n} \Big | \sum_{i=0}^{k-1} (f \circ T^i - \nu (f) -Z_i) \Big|= o(n^{1/p}L(n)) \quad
  \text{almost surely,}
\end{equation}
for $p\in ]2,4]$ and $L$ an explicit slowly varying function. For different classes of piecewise expanding maps $T$ of $[0,1]$, almost sure invariance principles with good remainder estimates have been established by Melbourne and Nicol (\cite{MelbourneNicol1}, \cite{MelbourneNicol2}) for
H\"older observables, and by Merlev\`ede and Rio \cite{MR} under rather mild
integrability assumptions. For instance, for uniformly expanding
maps as defined in  Definition \ref{UE}, Merlev\`ede and Rio \cite{MR} obtained in \eqref{ASIPp} a rate of order $n^{1/3} ( \log n)^{1/2} ( \log \log n)^{(1+ \varepsilon )/3}$ for any $\varepsilon >0$ and a class of
observables $f$ in ${\mathbb L}^r(\nu)$ where $r>3$.
 Because of the intrinsic time-orientation  of the non-invertible dynamical 
 systems studied in the papers above mentioned, the almost sure invariance principle cannot be obtained  directly  by mean of a martingale approximation as done for instance in \cite{FMT}.  In \cite{MR}, the approximating Brownian motion is constructed with the help of the conditional quantile transform combined with a coupling method based, roughly speaking, on a conditional version of the Kantorovitch-Rubinstein theorem. The approach used in \cite{MelbourneNicol1} or in \cite{MelbourneNicol2} when $d=1$ (see their appendix A) is based on the  seminal paper by  Philipp and Stout \cite{PS} which proves ASIP by combining martingale approximation with blocking arguments and the Strassen-Skorohod embedding. With this approach, they cannot reach better rates than $O(n^{\varepsilon + 3/8})$ in \eqref{ASIPp} for any $\varepsilon >0$ (see Remark 1.7 in \cite{MelbourneNicol2}). Notice that recently Gou\"ezel \cite{Go} obtained a very general result concerning the rates in the ASIP for \emph{vector-valued} observables of the iterates of dynamical systems by mean of spectral methods. In the particular case of expanding maps as defined in Definition \ref{UE}, his result gives the rate $o(n^{1/4 + \varepsilon})$ for any $\varepsilon >0$ for some bounded vector-valued observables.

On an other hand if we consider the strong approximation principle of the partial sums of real-valued functions of the Markov chain associated to the dynamical system, better rates can be reached. Let denote by $K$ the
Perron-Frobenius operator of $T$ with respect to
$\nu$. Recall that for any bounded measurable functions $f$ and $ g$,
\beq \label{defPF}
\nu(f \cdot  g\circ T)=\nu(K(f) g) \, .
\eeq
Let $(Y_i)_{i \geq 0}$ be a stationary Markov chain with invariant
measure $\nu$ and transition Kernel $K$. The sequence $(Y_i)_{i \geq 0}$  corresponds to the iteration of the inverse branches of
$T$. Combining a suitable martingale approximation with a sharp control of the approximation error in ${\mathbb L}^p$ for $p\geq 2$, and the Skorohod-Strassen embedding theorem as done in Shao \cite[Theorem 2.1]{shao93}, Dedecker, Merlev\`ede and Doukhan \cite{DDM}  obtained projective conditions under which the strong approximation holds with an explicit rate, as
 in \eqref{ASIPp}, for the partial sums associated to a stationary sequence of real valued random variables. In the Markov chain setting, the conditions involved in their corollaries 2.1 and 2.2 can be rewritten with the help of the transition kernel $K$. One can then wonder whether under the same conditions, the strong approximation principle also holds with the same rate for the partial sums associated to some observables of the iterates of the dynamical system. Indeed, it is well known (see for instance Lemma XI.3 in Hennion and Herv\'e \cite{HH})
that on the probability space $([0, 1], \nu)$, the random
vector $(T, T^2, \ldots , T^n)$ is
distributed as $(Y_n,Y_{n-1}, \ldots, Y_1)$. Hence any information on the law of
$\sum_{i=1}^n (f \circ T^i-\nu(f))$ can be
obtained by studying  the law of  $\sum_{i=1}^n (f(X_i)-\nu(f))$.
However, the reverse time property cannot be used to transfer directly the almost sure results for $\sum_{i=1}^n (f(Y_i)-\nu(f))$ to the sum $\sum_{i=1}^n (f \circ T^i-\nu(f))$. To prove results on the strong approximation principle for the partial sums of functions of the iterates of the dynamical system, it is then needed to work on the dynamical system itself. Therefore to prove that $\sum_{i=1}^n (f \circ T^i-\nu(f))$ satisfies a strong approximation principle with the same rates than the one reached for the   associated Markov chain, we shall approximate the partial sum associated to real-valued observables of the iterates of the dynamical system by a sum of reverse martingale differences $M^*_n$. To be more precise,  we shall approximate $S_n(f)$ by $M^*_n=\sum_{k=1}^n d_k^*$ where $(d_k^*)_{k \in {\mathbb N}}$ is a sequence of real valued random variables that is measurable with respect to a non-increasing sequence of $\sigma$-algebras say $({\mathcal G}_k)_{k \in {\mathbb N}}$ and such that $\E (d_k^* | {\mathcal G}_{k+1})=0$ almost surely. An analogue of Theorem 2.1 in Shao \cite{shao93} but in the context of the partial sums associated  to a sequence of reverse martingale differences, is then needed. Up to our knowledge, this  version for reverse martingale differences sequences does not exist in the literature. Starting from a reverse martingale version of the Skorohod-Strassen embedding theorem as done in Scott and Huggins \cite{HS}, we shall prove it in Section \ref{sectionreverseres} (see our Theorem \ref{shaoresult}).

Our paper is organized as follows. In Section \ref{sectionreverseres}, we state some results concerning the almost sure invariance principle with rates for sums associated to a sequence of reverse martingale differences. The proofs of these results are postponed to Section \ref{sectproofmartingale}. In Section \ref{sectionASIPUEM}, we obtain rates of convergence in the ASIP for the sum of a large class of functions (non necessarily bounded) of the iterates of some classes of uniformly expanding maps,  as well as for the partial sum associated to the corresponding Markov chain. A part of these results coming from a more general result on $\phi$-dependent sequences, we state and prove in Section \ref{sectiongeneralasip} an ASIP with explicit rates for functions of a stationary sequence satisfying a mild $\phi$-dependent condition. Section \ref{sectionproofsectionASIPUEM} is devoted to the proofs of the results on expanding maps stated in Section \ref{sectionASIPUEM}.

In this paper, we shall  use sometimes the notation $a_{n}\ll b_{n}$ to mean that there exists a numerical constant $C$
not depending on $n$ such that $a_{n}\leq C b_{n}$, for all positive integers $n$.

\section{ASIP with rates for sums of
differences of reverse martingales} \label{sectionreverseres}

\setcounter{equation}{0}

The next two results can be viewed as  reverse martingales analogues to a result of Shao \cite[Theorem 2.1]{shao93}. The proof of the next result follows from the Skorohod embedding of reverse martingales in Brownian motion as obtained in Scott and Huggins \cite{HS} together with an estimate for Brownian motions given in Hanson and Russo \cite[Theorem 3.2A]{HR}. The proofs of all the results presented in this section are postponed to Section \ref{sectproofmartingale}.

\medskip

In this section we consider real random variables defined on a probability
space $(\Omega,{\mathcal{A}},{\mathbb{P}})$.

\begin{prop}\label{huggins}
Let $(\xi_n)_{n \in {\mathbb N}}$ be a sequence of square integrable variables adapted
to a non-increasing filtration $(\G_n)_{n \in {\mathbb N}}$. Assume that $\E(\xi_n|\G_{n+1})=0$ a.s.
and $\delta_n^2:=\sum_{k\ge n} \E(\xi_k^2)<\infty$.
Then $V_n^2:=\sum_{k\ge n}
\E(\xi_k^2|\G_{k+1})$ is well defined a.s. and in ${\mathbb L}^2$ as well as
$R_n=\sum_{k\ge n} \xi_k$. Let $(\alpha_n)_{n \in {\mathbb N}}$ be a sequence of non-increasing
positive numbers with $\alpha_n =O(\delta_n^2)$ and $\alpha_n/\delta_n^4\to
\infty$. Assume that
\begin{gather}
\label{cond1} V_n^2-\delta_n^2=o(\alpha_n) \qquad \mbox{$\P$-a.s.} \\
\label{cond2}\sum_{n\ge 1} \alpha_n^{-\nu} \E(|\xi_n|^{2\nu}
) \quad \mbox{ for some $1\le \nu \le 2$} \, .
\end{gather}
Then, enlarging our probability space if necessary it is possible to
find a standard Brownian motion $(B_t)_{t\ge 0}$, such that
\begin{equation}\label{raterenv}
R_n-B_{\delta_n^2}=o \Big ( \big (\alpha_n\big(|\log (\delta_n^2/\alpha_n)|+\log \log (\alpha_n^{-1})\big) \big )^{1/2}\Big ) \qquad \mbox{$\P$-a.s.}
\end{equation}
\end{prop}
\begin{rmk} \label{rmkO}It follows from the proof that if \eqref{cond1} holds with "big $O$'' instead of "little $o$'' then \eqref{raterenv} holds with the same change.
\end{rmk}

\medskip

Now, we give  a result for partial sums associated to a sequence of reverse martingale differences 
rather than for tail series. It may be viewed as the analogue of Theorem 2.1 
in Shao \cite{shao93} but in the context of reverse martingale differences.

\begin{thm} \label{shaoresult}
Let $(X_n)_{n \in {\mathbb N}}$ be a  sequence of square integrable random variables adapted to a non-increasing filtration  $(\G_n)_{n \in {\mathbb N}}$. Assume that $\E(X_n|\G_{n+1})=0$ a.s., that $\sigma_n^2:=\sum_{k=1}^n \E(X_k^2)\to \infty$ and that
$\sup_n \E(X_n^2)<\infty$.
Let $(a_n)_{n \in {\mathbb N}}$ be a non-decreasing sequence of positive numbers such that
$(a_n/\sigma_n^2)_{n \in {\mathbb N}}$ is non-increasing and $(a_n/\sigma_n)_{n \in {\mathbb N}}$ is non-decreasing.
Assume that
\begin{gather}\label{cond1'}
\sum_{k=1}^n (\E(X_k^2|\G_{k+1})-\E(X_k^2))=o(a_n) \qquad \mbox{$\P$-a.s.} \\
\label{cond2'}\sum_{n\ge 1} a_n^{-\nu} \E(|X_n|^{2\nu})<\infty \quad \mbox{ for some $1\le \nu \le 2$} \, .
\end{gather}
Then, enlarging our probability space if necessary it is possible to
find a sequence $(Z_k)_{k \geq 1}$ of independent centered  Gaussian variables with
$\E(Z_k^2)=\E(X_k^2)$ such that
\begin{equation}\label{rate}
\sup_{1 \leq k \leq n } \Big | \sum_{i=1}^k X_i - \sum_{i=1}^k Z_i \Big | = o \big(\big ( a_n(|\log (\sigma_n^2/a_n)|+\log \log a_n)
\big)^{1/2}\big)
\qquad \mbox{$\P$-a.s.}
\end{equation}
\end{thm}
\begin{rmk}An inspection of the proof allows to weaken slightly some of the assumptions as follows. Assume that $\E(X_n^2)=O(\sigma_n^{2s})$
for some $0\le s <1$ instead of boundedness and assume that there exists $C>1$,
such that for every $n\ge 1$,  $\sup_{k\ge n} (a_k/\sigma_k^2)\le C a_n/
\sigma_n^2$, and $\inf_{k\ge n} (a_k/\sigma_k)\ge a_n/(C\sigma_n)$, instead of
the corresponding monotonicity assumptions.
\end{rmk}

\medskip

We derive now the functional LIL for the partial sums associated to a stationary  and ergodic sequence of reverse martingale differences. In the next corollaries, we make use of $\theta:\Omega\mapsto\Omega$  a measurable transformation preserving the probability ${\mathbb{P}}$.
\begin{cor} \label{corFLIL}
Let $X_0$ in ${\mathbb L}^2$ and for $n \geq 1$, $X_n = X_0 \circ \theta^n$.  Assume that $\theta$ is ergodic. Let
$(\G_n)_{n \geq 0}$ be a non-increasing filtration to which $(X_n)_{n\geq 0}$ is adapted and
such that $\E(X_n|\G_{n+1})=0$ a.s.  Enlarging our probability space if necessary it is possible to
find a sequence $(Z_k)_{k \geq 1}$ of independent centered  Gaussian variables with
$\E(Z_k^2)=\E(X_1^2)$ such that
\begin{equation}\label{FLIL}
\sup_{1 \leq k \leq n }\Big | \sum_{i=1}^k X_i - \sum_{i=1}^k Z_i \Big | = o \big ( \sqrt {n\log \log n} \big )
\qquad \mbox{$\P$-a.s.}
\end{equation}
\end{cor}
\begin{rmk}The Strassen functional law of the iterated logarithm (FLIL) follows
from the corollary. Notice that a semi-FLIL has been
proved by Wu \cite{Wu}.
\end{rmk}
\medskip

We now give rate results in the strong invariance principle for the partial sums associated to a stationary   sequence of reverse martingale differences.

\begin{cor} \label{corrateASIP}
Let $2<p <4$. Let $X_0$ be in ${\mathbb L}^p$ and for $n \geq 1$, $X_n = X_0 \circ \theta^n$.   Let
$(\G_n)_{n \geq 0}$ be a non-increasing filtration to which $(X_n)_{n\geq 0}$ is adapted and
such that $\E(X_n|\G_{n+1})=0$ a.s. Let $b(\cdot)$ be a positive non-decreasing slowly
varying function, such that $x \mapsto x^{2/p-1} b(x)$ is non-increasing. Assume that
\begin{equation} \label{condcarrecondi}
\sum_{k=1}^n (\E(X_k^2|\G_{k+1})-\E(X_k^2)) = o(n^{2/p}b(n)) \qquad \as
\end{equation}
Enlarging our probability space if necessary it is possible to
find a sequence $(Z_k)_{k \geq 1}$ of independent centered  Gaussian variables with
$\E(Z_k^2)=\E(X_1^2)$ such that
\begin{equation}\label{rateFLIL}
\sup_{1 \leq k \leq n }\Big | \sum_{i=1}^k X_i - \sum_{i=1}^k Z_i \Big |= o \big ( n^{1/p}\sqrt {b(n)\log n} \big )
\qquad \mbox{$\P$-a.s.}
\end{equation}
\end{cor}
\begin{cor} \label{corrateASIPp=4}
Let ${\mathcal
G}_0$ be a sub-$\sigma$-algebra of $\mathcal{A}$ satisfying
$  \theta^{-1}({\mathcal G}_0) \subseteq {\mathcal G}_0$. Let $X_0$ be a ${\mathcal
G}_0$-measurable random variable in ${\mathbb L}^4$. For $n \geq 1$,  let $X_n = X_0 \circ \theta^n$ and $\G_n = \theta^{-n} ({\mathcal G}_0)$. Assume that
$\theta$ is ergodic and that $\E(X_n|\G_{n+1})=0$ a.s.
Assume that
\begin{equation*}
\sum_{k=1}^n (\E(X_k^2|\G_{k+1})-\E(X_k^2)) = O(( n \log \log n )^{1/2}) \qquad \as
\end{equation*}
Enlarging our probability space if necessary it is possible to
find a sequence $(Z_k)_{k \geq 1}$ of independent centered  Gaussian variables with
$\E(Z_k^2)=\E(X_1^2)$ such that
\begin{equation}\label{rateFLILp=4}
\sup_{1 \leq k \leq n }\Big | \sum_{i=1}^k X_i - \sum_{i=1}^k Z_i \Big |= O \big ( n^{1/4} (\log n)^{1/2} (\log \log n)^{1/4}\big )
\qquad \mbox{$\P$-a.s.}
\end{equation}
\end{cor}
\begin{rmk} \label{remarkdirectsense}
If we consider usual martingale differences, so more precisely if $(X_n)_{n \in {\mathbb N}}$ is  a stationary sequence
in ${\mathbb L}^p$ for a real $p \in ]2,4]$ adapted to a  non-decreasing filtration $(\F_n)_{n \in {\mathbb N}}$ (assume moreover that $\F_n=\theta^{-n}\F_0$ when 
$p=4$) and such that $\E(X_n|\F_{n-1})=0$ a.s., the usual Skorohod-Strassen embedding theorem gives  similar results to Corollaries \ref{corrateASIP} and \ref{corrateASIPp=4}. Indeed starting from Theorem 2.1 in Shao \cite{shao93} when $p\in ]2,4[$ or from its proof when $p=4$ (notice that by stationarity -of the process and the filtration-, the stopping times used to construct the approximating Brownian motion can be chosen to be stationary), and using the arguments developed to get Corollaries \ref{corrateASIP} and \ref{corrateASIPp=4}, we infer that if there exists a positive non-decreasing slowly varying function $b(\cdot)$, such that
\begin{equation*}
\sum_{k=1}^n (\E(X_k^2|\F_{k-1})-\E(X_k^2)) = o(n^{2/p}b(n)) \qquad \as \, \,\mbox{when $p\in ]2,4]$},
\end{equation*}
or if, $(X_n)_{n\in \mathbb{N}}$ is ergodic and
\begin{equation*}
\sum_{k=1}^n (\E(X_k^2|\F_{k-1})-\E(X_k^2)) = O(( n \log \log n )^{1/2}) \qquad \as \, \, \mbox{when $p=4$},
\end{equation*}
then, enlarging our probability space if necessary it is possible to
find a sequence $(Z_k)_{k \geq 1}$ of independent centered  Gaussian variables with
$\E(Z_k^2)=\E(X_1^2)$ such that the rates of approximation in the strong invariance principle are the same as those given by \eqref{rateFLIL} and \eqref{rateFLILp=4}.
\end{rmk}

\section{ASIP with rates for uniformly expanding maps} \label{sectionASIPUEM}
 \setcounter{equation}{0}
Several classes of uniformly expanding maps of the interval are
considered in the literature. In Theorem \ref{ASmap1} below, we shall consider the following definition:
\begin{defn}\label{UE}
A map $T:[0,1] \to [0,1]$ is uniformly expanding if it belongs to the
class ${\mathcal C}$ defined in Broise \cite[Section 2.1 page 11]{Br}
and if
\begin{enumerate}
\item There is  a unique absolutely continuous invariant probability
measure $\nu$, whose density $h$ is such that $\frac{1}{h}{\bf
1}_{h>0}$ has bounded variation.
\item The system $(T, \nu)$ is mixing in the ergodic theoretic sense.
\end{enumerate}
\end{defn}
We refer also to Definition 1.1 in \cite{DeGoMe} for a more complete description. Some well known examples of maps satisfying the
above conditions are:
\begin{enumerate}
\item $T(x)= \beta x -[\beta x]$ for $\beta> 1$. These maps are called
$\beta$-transformations.
\item  $I$ is the finite union of disjoint intervals $(I_k)_{1 \leq k
\leq n}$, and $T(x)=a_kx +b_k$ on $I_k$, with $|a_k|>1$.
\item $T(x)=a(x^{-1}-1)-[a(x^{-1}-1)]$ for some $a>0$. For $a=1$, this
transformation is known as the Gauss map.
\end{enumerate}

\smallskip

Our aim is to obtain explicit rates in the strong invariance principle for $\sum_{i=0}^{n-1} ( f \circ T^i - \nu (f))$ and $\sum_{i=1}^n ( f (Y_i) - \nu (f))$ ($(Y_i)_{i\geq 0}$ corresponding to the iterations of the inverse branches of $T$) when $T$ is a dynamical system defined in Definition \ref{UE}, and $f$ belongs to a large class of functions non necessarily bounded. Such
classes are described in the following definition.
\begin{defn}
If $\mu$ is a probability measure on $\mathbb R$ and $p \in [2,
\infty)$, $M \in (0, \infty)$, let $\Mon_p(M,\mu)$ denote the set of
functions $f:\R\to \R$ which are monotonic on some interval and null
elsewhere and such that $\mu(|f|^p)\leq M^p$. Let $\Monm_p(M,\mu)$ be
the closure in ${\mathbb L}^p(\mu)$ of the set of functions which can
be written as $\sum_{\ell=1}^L a_\ell f_\ell$, where $\sum_{\ell=1}^L|a_\ell| \leq 1$ and $f_\ell\in \Mon_p(M, \mu)$.
\end{defn}
\begin{rmk}
In previous papers, see for instance \cite{DeGoMe}, the closure in $\mathbb{L}^1(\mu)$ was used in the definition above.
It turns out that both definitions coincide. Indeed, a sequence bounded in $\mathbb{L}^p(\mu)$ and converging in
$\LL^1(\mu)$, converges for the weak topology in $\LL^p(\mu)$. To conclude, recall that, by the Hahn-Banach theorem,
in any Banach space, the weak closure of a convex set is equal to its strong closure.
\end{rmk}

Our main theorem follows. For uniformly expanding maps as defined in Definition \ref{UE}, it involves
an ${\mathbb L}^p$-integrability condition of the observables.

\begin{thm} \label{ASmap1} Let $T$ be a uniformly expanding  map as defined in Definition \ref{UE},
with absolutely continuous invariant measure $\nu$. Let $p\in ]2,4]$. Then, for any
$M>0$ and any $f \in \Monm_p(M,\nu)$, the series
\begin{equation}\label{var}
\sigma^2=\sigma^2(f)= \nu((f-\nu(f))^2)+ 2 \sum_{k>0}
\nu ((f-\nu(f))f\circ T^k)
\end{equation}
converges absolutely to some non negative number.
\begin{enumerate}
\item Let $(Y_i)_{i \geq 1}$ be a stationary Markov chain with
transition kernel $K$ defined by \eqref{defPF} and invariant measure $\nu$. Enlarging the underlying probability space if necessary, there exists a
    sequence $(Z_i)_{i \geq 0}$ of iid Gaussian random
    variables with mean zero and variance $\sigma^2$ defined by
    \eqref{var}, such that
\begin{equation*}
\sup_{1 \leq k \leq n}\Big|\sum_{i=1}^{k} (f (Y_i)
  -\nu(f))- Z_i) \Big| = \left\{ \begin{array}{ll} o (n^{1/p} ( \log n)^{1 - 2/p} ) \mbox{ a.s. } & \mbox{ if $p \in ]2,4[$}  \\
O(n^{1/4} ( \log n)^{1/2} ( \log \log n)^{1/4}) \mbox{ a.s. } & \mbox{ if $p=4$}.
\end{array}
\right.
\end{equation*}
\item Enlarging the probability space $([0, 1], \nu)$ if necessary,  there exists a
    sequence $(Z^*_i)_{i \geq 0}$ of iid Gaussian random
    variables with mean zero and variance $\sigma^2$ defined by
    \eqref{var}, such that
\begin{equation*}
\sup_{1 \leq k \leq n}\Big|\sum_{i=0}^{k-1} (f \circ T^i
  -\nu(f))- Z^*_i) \Big| = \left\{ \begin{array}{ll} o (n^{1/p} ( \log n)^{1 - 2/p} ) \mbox{ a.s. } & \mbox{ if $p \in ]2,4[$}  \\
O(n^{1/4} ( \log n)^{1/2} ( \log \log n)^{1/4}) \mbox{ a.s. } & \mbox{ if $p=4$}.
\end{array}
\right.
\end{equation*}
\end{enumerate}
\end{thm}
As we shall see in the proof of Theorem \ref{ASmap1}, Item 1 will follow from a more general result on $\phi$-dependent sequences. This general result is presented in a separate section (see Section \ref{sectiongeneralasip}) since it has an interest in itself. Item 2 is obtained by considering an approximation by partial sums associated to a sequence
 (which is non stationary when $p \in ]2,4[$) of reverse martingale differences and by applying the results of Section \ref{sectionreverseres}. As we shall see, the reverse time property mentioned in the introduction allows to make links  with estimates obtained to prove Item 1. Notice that when $p=2$, the strong invariance principle (therefore with the rate $o(n^{1/2} (\log\log n)^{1/2} )$) has been proved recently in \cite{DeGoMe} (see Item 3 of their Theorem 1.5). Notice that the proof of the strong invariance principle obtained in \cite{DeGoMe} combines different approximation arguments. The observable is approximated by a function with better integrability
properties, the almost sure invariance principle with rates obtained in \cite{MR} is then used for this new function, and a bounded law of
the iterated logarithm is proved to make it possible to pass the results from the
better function to the original function. Let us emphasize that the strong invariance principle obtained in \cite{DeGoMe} could be also proved by using a direct approximation by partial sums associated to a sequence of reverse martingale differences as we do in the proof of Theorem \ref{ASmap1} above, and by using our Corollary \ref{corFLIL}.

\medskip

Let us mention that when $p=4$, Theorem \ref{ASmap1}  is a consequence of  the following more
general result.

\begin{thm} \label{ASmapMC} Let $T$ be a map from $[0,1]$ to $[0,1]$ preserving a probability $\nu$  on $[0,1]$, and let $K$ be defined by \eqref{defPF}. Let $f$ be a measurable function such that $\nu (f^4) < \infty$. Assume that there exists $\gamma \in ]0,1]$ such that
\beq \label{cond1ASmapMC}
 \sum_{n>0}   (\log n)^3 n^{\frac{1}{\gamma} + \frac{1}{2}}  \| K^n (f) - \nu (f) \|^{2}_{4,\nu}  < \infty \, ,
\eeq
and
\beq \label{cond2ASmapMC}
 \sum_{n>0}  (\log n)^3 n^{2 \gamma}  \sup_{i \geq j \geq n }\| K^j (f K^{i-j}(f)) - \nu (f K^{i-j}(f) ) \|^{2}_{2, \nu } < \infty \, .
\eeq Then, the series $\sigma^2$ defined by \eqref{var} converges absolutely and both items of Theorem \ref{ASmap1} holds with rate $O (
 n^{1/4}  ( \log n)^{1/2} ( \log \log n)^{1/4} ) $.
\end{thm}
\begin{rmk}
We would like to emphasize that the strategy of proof of Item 2 for both Theorems \ref{ASmap1} and \ref{ASmapMC} is to approximate $S_n(f)=\sum_{i=0}^{n-1} (f \circ T^i
  -\nu(f))$ by a partial sum associated to a sequence of reverse martingale differences (possibly non stationary), let say $M_n^*=\sum_{i=0}^{n-1} d_i^*$,  to apply our strong approximation results given in Section \ref{sectionreverseres} and to have a nice control of the approximation error between 
 $S_n(f)$ and $M_n^*$.
 A careful analysis of the proof of Theorem 2.4 and its Corollary 2.7 in \cite{DDM} together with the arguments developed in the proof of Item 2 of our Theorems \ref{ASmap1} and \ref{ASmapMC}, show that our strategy of proof also gives new results when we consider classes of expanding maps with a neutral fixed point at zero such as the generalized Pomeau-Manneville (GPM) map as defined in Definition 1.1 of \cite{DeGoMe10}. Therefore we infer that if $T$ is a GPM map with parameter $\gamma \in (0,1)$ and $f$ is a function satisfying Condition (3.22) in \cite{DDM}, then the conclusion of Corollary 3.18 in \cite{DDM} also holds for $S_n(f)$. In particular, if $f$ is a bounded variation function and $T$ is a GPM map with parameter $\gamma \in (0,\delta^{-1}]$ where $\delta= p+1 -2/p$ and $p \in ]2, 4]$, then  $S_n(f)$ satisfies an almost sure invariance principle with rate $o(n^{1/p} \log n)$.

\end{rmk}
For a map $T$ from $[0,1]$ to $[0,1]$ preserving a probability $\nu$  on $[0,1]$ and $(Y_i)_{i \geq 0}$ its associated Markov chain with invariant
measure $\nu$ and transition Kernel $K$ defined by \eqref{defPF}, Theorem \ref{ASmapMC} also allows to obtain rates in the strong invariance principle for  $\sum_{i=1}^n (f \circ T^i-\nu(f))$ (or for the partial sum of its associated Markov chain) when $f$ has a modulus of continuity that is dominated by a concave and non-decreasing function and the condition \eqref{cond1PFO} below is satisfied.

 For any integer $k$, we denote by $Q_k$ the operator defined as follows: $\E ( g ( Y_0,  Y_k) | Y_0 = x) = Q_k (g) (x)$. Let $\Lambda_1({\mathbb R})$ be the set of functions from ${\mathbb R}$ to ${\mathbb R}$ that are $1$-Lipschitz and let $\Lambda_1({\mathbb R}^2)$ be the set of functions $h$
from ${\mathbb R}^2$ to ${\mathbb R}$ such that
$$
  |h(x_1, y_1)-h(x_2, y_2)|\leq \frac 12 |x_1-x_2|+ \frac 12 |y_1-y_2| \, .
$$
Denote by $\|\cdot \|_{\infty, \mu}$ the essential supremum
norm with respect to $\nu$. Assume that there exist $C>0$ and $\rho \in ]0, 1[$ such that for  any $(i,j) \in {\mathbb N}^2$,
\beq \label{cond1PFO}
\sup_{g \in \Lambda_1({\mathbb R})} \Vert K^i (g)
- \nu (g) \Vert_{\infty, \nu}\leq C \rho^i  \ \text{ and } \
\sup_{j \geq 0} \sup_{h \in \Lambda_1({\mathbb R}^2)} \Vert K^i \circ Q_j  (h)
- \nu \big ( Q_j  (h) \big ) \Vert_{\infty, \nu}\leq C \rho^i \, .
\eeq
\begin{thm} \label{ASmapMCconcave}
Let $T$ be a  map from $[0,1]$ to $[0,1]$ preserving a probability $\nu$  on $[0,1]$ and  let $(Y_i)_{i \geq 0}$ be a stationary Markov chain with invariant
measure $\nu$ and transition Kernel $K$ defined by \eqref{defPF}. Assume that condition \eqref{cond1PFO} is satisfied.  Let $f$ from ${\mathbb R}$ to ${\mathbb
R}$ such that $ | f(x) - f(y) | \leq c( |x-y|) $, for some concave
and non-decreasing function $c$ satisfying
\begin{equation}\label{intc}
\int_0^1 |\log t |^{2/{\sqrt 3}} | \log \log (2 t^{-1}) |^3 \frac{c(t)}{t } dt < \infty \, .
\end{equation}
Then the conclusion of Theorem \ref{ASmapMC} holds.
\end{thm}
Note that (\ref{intc}) holds if $c(t)\leq D|\log(t)|^{-\gamma}$ for
some $D>0$ and some $\gamma>1 + 2/{\sqrt 3}$. Therefore Theorem \ref{ASmapMCconcave} applies to the  functions from $[0,1]$ to ${\mathbb R}$ which
are $\alpha$-H\"older for some $\alpha \in ]0, 1]$.

\medskip

Dedecker and Prieur \cite[Section 7.2, Example 4.4]{DP07} have shown that condition \eqref{cond1PFO} is satisfied for a
large class of uniformly expanding maps such as those considered in \cite{CMS}. The conditions imposed on the class of the expanding maps they consider are slightly more restrictive than those considered in Definition \ref{UE}. In particular they are defined by mean of finite partitions of $[0,1]$.

\section{A general ASIP result for a class of weakly dependent sequences} \label{sectiongeneralasip}

 \setcounter{equation}{0}

Let $(\Omega ,\mathcal{A}, \p)$ be a probability space, and let
$\theta :\Omega \mapsto \Omega $ be a bijective bimeasurable
transformation preserving the probability ${\p}$. Let ${\mathcal
F}_0$ be a sub-$\sigma$-algebra of $\mathcal{A}$ satisfying
${\mathcal F}_0 \subseteq \theta^{-1}({\mathcal F}_0)$. Define the nondecreasing
filtration $({\cal F}_i)_{i \in {\mathbb Z}}$ by
${\cal F}_i =\theta^{-i}({\mathcal F}_0)$. We shall sometimes denote  by ${\mathbb E}_i$ the conditional expectation with respect to ${\mathcal F}_i$, and we shall set $P_i ( \cdot) = {\mathbb E}_i ( \cdot) - {\mathbb E}_{i-1}( \cdot) $.

\begin{defn}\label{defphi}
For any integrable random variable $X$, let us write
$X^{(0)}=X- \E(X)$.
For any random variable $Y=(Y_1, \cdots, Y_k)$ with values in
${\mathbb R}^k$ and any $\sigma$-algebra $\F$, let
\[
\phi(\F, Y)= \sup_{(x_1, \ldots , x_k) \in {\mathbb R}^k}
\Big  \| \E \Big(\prod_{j=1}^k (\I_{Y_j \leq x_j})^{(0)} \Big | \F \Big)^{(0)}
\Big \|_\infty.
\]
For a sequence ${\bf Y}=(Y_i)_{i \in {\mathbb Z}}$, where $Y_i=Y_0
\circ \theta^i$ and $Y_0$ is an $\F_0$-measurable and real-valued
random variable, let
\begin{equation*}
\phi_{k, {\bf Y}}(n) = \max_{1 \leq l \leq
k} \ \sup_{ n\leq i_1\leq \ldots \leq i_l} \phi(\F_0,
(Y_{i_1}, \ldots, Y_{i_l})) \, .
\end{equation*}
\end{defn}
\begin{thm}\label{asipgenephi}
Let  $X_i = f(Y_i) - \E ( f(Y_i))$, where $Y_i=Y_0 \circ \theta^i$
and $Y_0$ is an $\F_0$-measurable random variable. Let $P_{Y_0}$ be the distribution of $Y_0$ and  $p \in ]2,4]$. Assume that $f$ belongs to $\Monm_p(M,P_{Y_0})$ for some $M>0$ and that
\begin{equation}\label{condDDM}
\sum_{k \geq 1} k^{1/\sqrt 3 -1/2} \phi^{1/2}_{2, {\bf Y}}(k) < \infty \  \text{ if $p \in ]2 , 4[ \, $ and } \ \sum_{k \geq 2}  (\log k)^3 k^{2/ {\sqrt 3} } \phi_{2, {\bf Y}}(k) < \infty \  \text{ if $p = 4$ }\, .
\end{equation}
Then, enlarging our probability space if necessary, there exists a  
    sequence $(Z_i)_{i \geq 0}$ of iid\ Gaussian random
    variables with mean zero and variance $\sigma^2$ defined by the absolutely converging series
\beq \label{varweakly}
\sigma^2 = \sum_{k \in {\mathbb Z}} {\rm Cov} ( X_0,X_k)
\eeq
such that
\begin{equation*}
\sup_{1 \leq k \leq n }\Big|\sum_{i=1}^{k} (X_i
- Z_i) \Big|= \left\{ \begin{array}{ll} o (n^{1/p} ( \log n)^{1 - 2/p} ) \mbox{ a.s. } & \mbox{ if $p \in ]2,4[$}  \\
O(n^{1/4} ( \log n)^{1/2} ( \log \log n)^{1/4}) \mbox{ a.s. } & \mbox{ if $p=4$}.
\end{array}
\right.
\end{equation*}
\end{thm}
Notice that an application of Corollary 2.1 in \cite{DDM} would  give a	 rate of convergence of order $o( n^{1/p} ( \log n)^{(t+1)/2})$ with $t>2/p$. As we shall see in the proof of the above theorem, the improvement of
the power in the logarithmic term is achieved via some truncation arguments.



The proof makes use of the following lemma (see Lemma 5.2 in Dedecker, Gou\"ezel and Merlev\`ede \cite{DeGoMe}).
\begin{lma}\label{covaphi}
Let ${\bf Y}=(Y_i)_{i \in {\mathbb Z}}$, where $Y_i=Y_0 \circ
\theta^i$ and $Y_0$ is an $\F_0$-measurable random variable. Let $f$
and $g$ be two functions from ${\mathbb R}$ to ${\mathbb R}$ which
are monotonic on some interval and null elsewhere. Let $p \in [1,
\infty]$. If $\|f(Y_0)\|_p < \infty$, then, for any positive integer
$k$,
\[
  \|{\mathbb E}(f(Y_k)|{\mathcal F}_0)-{\mathbb E}(f(Y_k))\|_p \leq 2
  (2\phi_{1, {\bf Y}}(k))^{(p-1)/p}\|f(Y_0)\|_p \, .
\]
If moreover $p\geq 2$ and $\|g(Y_0)\|_p < \infty$, then for any positive integers
 $i\geq j\geq k$,
\[
\|{\mathbb E}(f(Y_i)^{(0)}g(Y_j)^{(0)}|{\mathcal F}_0)-
{\mathbb E}(f(Y_i)^{(0)}g(Y_j)^{(0)})\|_{p/2}\leq 8
 (4\phi_{2, {\bf Y}}(k))^{(p-2)/p}\|f(Y_0)\|_{p}\|g(Y_0)\|_p \, .
\]
\end{lma}

\subsection{Proof of Theorem \ref{asipgenephi} for $2<p<4$}

Let  $f\in \Monm_p(M,P_{Y_0})$. We shall first prove that
$\sum_{\ell \geq 0} \Vert \E_0(X_{\ell} )\Vert_p< \infty$ provided that  $\sum_{k \geq 1}  \phi^{(p-1)/p}_{1, {\bf Y}}(k) < \infty$ (notice that this condition is implied by \eqref{condDDM}). This will entail that $d_0:=\sum_{\ell \ge 0}P_0(X_\ell)$ and $r_0:=\sum_{\ell\ge 1} \E_0(X_\ell)$  are well defined and that  we have the
martingale-coboundary decomposition
 $X_0=d_0+r_0\circ \theta^{-1}-r_0$. Then the series $\sigma^2= \sum_{k \in {\mathbb Z}} {\rm Cov} ( X_0,X_k)$  will converge absolutely and we will have  $ \lim_{n \rightarrow \infty }n^{-1} \E (S^2_n(f)) = \sigma^2$. \\Since $f\in \Monm_p(M,P_{Y_0})$, by definition, there exists a sequence
of functions
\beq \label{defifL}
f_L= \sum_{k=1}^L a_{k, L} f_{k, L} \, ,
\eeq with
$f_{k, L}$ belonging to $\Mon_p(M,P_{Y_0})$ and $\sum_{k=1}^L
|a_{k, L}| \leq 1$, such that $f_L$ converges in ${\mathbb
L}^p(P_{Y_0})$ to $f$. Hence, $\| \E_0(X_\ell)\|_p=\lim_{L\to \infty}\|\E_0 (  f_L(Y_\ell)-\E(f_L(Y_\ell))\|_p\le \liminf_{L\to \infty}
\sum_{k=1}^L|a_{k,L}|\, \|\E_0(f_{k,L}(Y_{\ell})-\E(f_{k,L}(Y_\ell))\|_p$. Next, by Lemma \ref{covaphi}, we have $\|\E_0(f_{k,L}(Y_{\ell})-\E(f_{k,L}(Y_\ell))\|_p\le
2  (2\phi_{1, {\bf Y}}(k))^{(p-1)/p}M$. This shows that   $\sum_{\ell \geq 0}
\Vert \E_0(X_{\ell} )\Vert_p< \infty$ as soon as  $\sum_{k \geq 1}  \phi^{(p-1)/p}_{1, {\bf Y}}(k) < \infty$.

\medskip

 We shall  define random variables by truncating the functions $f_{k,L}$ defined
in \eqref{defifL}. Let  $j$ be a positive integer and
\[
g_j(x) = x {\bf 1}_{|x| \leq c(j)} \text{ where $c(j)=2^{j/p} j^{-2/p}$} \, .
\]
Define then
\beq \label{deffLj}
f_{L,j}=\sum_{k=1}^L a_{k, L} \,  g_j \circ f_{k, L}\, .
\eeq
By assumption and by construction, for any integers $j,L$, $f_{L,j}$  has a variation bounded by $4c(j)$. Hence, by
Lemma 2.1 of \cite{DeGoMe}, $(f_{L,j})_L$ admits a subsequence, say $(f_{\varphi(L),j})_L$ converging in $\LL^1(P_{Y_0})$, hence
in  $\LL^p(P_{Y_0)})$, say to $\bar f_j$.
Then,   $f-\bar f_j$ is the limit in ${\mathbb L}^p(P_{Y_0})$ of
\[
f_{\varphi(L)}-f_{\varphi(L), j}=
\sum_{k=1}^{\varphi(L)} a_{k, \varphi(L)}
\widetilde g_j \circ f_{k, L} \text{ where $\widetilde g_j = x {\bf 1}_{|x| > c(j)}$}\, .
\]
Define then for any integer $\ell$,
\beq \label{defXj}
\bar X_{j,\ell} := \bar f_j (Y_\ell) - \E (\bar f_j (Y_\ell)) \text{ and } \widetilde X_{j,\ell} := (f- \bar f_j) (Y_\ell) - \E ( (f -\bar f_j )(Y_\ell)) \, ,
\eeq
and
\beq \label{defbartildeXjL}
\bar X_{j,L, \ell} := f_{\varphi(L), j}(Y_\ell)  - \E (f_{\varphi(L), j}(Y_{\ell})) \text{ and }  \widetilde X_{j,L, \ell}
:= (f_{\varphi(L)} -f_{\varphi(L), j})(Y_\ell)  - \E ((f_{\varphi(L)}-f_{\varphi(L), j})(Y_{\ell})) \, .
\eeq
We define also a sequence of martingale differences, $(\bar d_{j,\ell})_{\ell \geq 1}$,  with respect to the non-decreasing sequence of  $\sigma$-algebras $({\mathcal F}_\ell)_{\ell \geq 1}$, as follows:
\beq \label{defdj}
\bar d_{j,\ell} = \sum_{k \geq \ell}P_\ell\left(\bar X_{j,k}\right)   \, .
\eeq
Notice that by assumption \eqref{condDDM}, the series $\sum_{k \geq 0}P_0\left(\bar X_{j,k}\right)$ converges in ${\mathbb L}^{\infty }$ as shown by the following
claim whose proof is given later.
\begin{claim} \label{convlqdjk}
Let $j$ be fixed. Assume that $ \sum_{k \geq 1} \phi_{1, {\bf Y}}(k) < \infty$.
Then $\sum_{k \geq 0} \Vert \E_0 \left(\bar X_{j,k}\right) \Vert_{\infty}< \infty $, and the sequence $(d_{j,\ell})_{\ell \geq 1}$ defined by \eqref{defdj} forms a stationary sequence of martingale differences in ${\mathbb L}^{\infty}$ with respect to the non-decreasing sequence of  $\sigma$-algebras $({\mathcal F}_\ell )_{\ell \geq 1}$.
\end{claim}
We define now some non stationary sequences $(\bar X_{\ell})_{\ell \geq 1}$ and $(\bar d_{\ell})_{\ell \geq 1}$ as follows:
\beq \label{def*1}
\bar d_1:=\bar d_{1,1}  \, , \  \bar X_1:=\bar X_{1,1} \, ,
\eeq
and, for every $j\ge 0$ and every $\ell\in\{2^{j}+1,...,2^{j+1}\}$,
\beq \label{def*} \bar d_\ell:=\bar d_{j,\ell}  \, , \
        \bar X_\ell:=\bar X_{j,\ell} \, .
\eeq
For every positive integer $n$, we then define
$$
\bar M_n(f):=\sum_{\ell=1}^n \bar d_\ell\ \ \mbox{and}\ \ \bar S_n(f):=\sum_{\ell=1}^n \bar X_{\ell} \, .$$
With these notations, the following decomposition is valid: for any positive integer $k$,
$$
S_k(f) = (S_k(f)  - \bar S_k(f))  + (\bar S_k(f)  - \bar{M}_k(f)) + \bar M_k(f) \, .
$$
Therefore, the theorem will follow if we can prove that
\beq \label{conv1thgene}
\sup_{1 \leq k \leq n } \big |S_k(f) - \bar S_k(f)\big | = o (n^{1/p} ( \log n)^{1 - 2/p} ) \text{ almost surely,}
\eeq
\beq \label{conv2thgene}
\sup_{1 \leq k \leq n } \big |\bar S_k(f)- \bar M_k(f)\big | = o (n^{1/p} ( \log n)^{1 - 2/p} ) \text{ almost surely,}
\eeq
and if, enlarging our probability space if necessary, there exists  a
    sequence $(Z_i)_{i \geq 0}$ of iid\ Gaussian random
    variables with mean zero and variance $\sigma^2$
such that
\begin{equation}\label{conv3thgene}
\sup_{1 \leq k \leq n }\Big|\sum_{i=1}^{k} (\bar d_i
- Z_i) \Big| =o (n^{1/p} ( \log n)^{1 - 2/p} )\,
  \text{almost surely.}
\end{equation}

\noindent{\sc Proof of \eqref{conv1thgene}.} For any non negative integer $j$, let
\begin{eqnarray*} \label{0dec}
\widetilde D_j := \sup_{ 1 \leq k \leq 2^{j}} |  \sum_{\ell= 2^{j}+1}^{k+2^{j}}(X_\ell-\bar X_\ell)| = \sup_{ 1 \leq k \leq 2^{j}} |  \sum_{\ell= 2^{j}+1}^{k+2^{j}}\widetilde X_{j,\ell}|  \, ,
\end{eqnarray*}
where $\widetilde X_{j,\ell}$ is defined by \eqref{defXj}. Let $N \in {\mathbb N}^*$ and let $k \in ]1, 2^N]$. We first notice that $\widetilde D_j \geq |  \sum_{\ell= 2^{j}+1}^{2^{j+1}}\widetilde X_{j,\ell}|$, so if $K$ is the  integer such that $2^{K-1} < k \leq  2^{K}$, then
 $$\big  | S_k(f) - \bar S_k(f) \big | \leq |X_1-\bar X_1|+\sum_{j=0}^{K-1} \widetilde D_j\, .$$
Consequently since $K \leq N$,
\begin{equation} \label{1dec} \sup_{1\le k\le 2^N}|S_k(f) - \bar S_k(f)|\le |X_1-\bar X_1|+\sum_{j=0}^{N-1}\widetilde D_j \, .
    \end{equation}
Therefore, \eqref{conv1thgene} will follow if we can prove that
$$
\sup_{1\le k\le 2^{j}}
    \Big |\sum_{\ell= 2^{j}+1}^{k+2^{j}}\widetilde X_{j,\ell}\Big | =o \big (  j^{1 - 2/p}  \, 2^{j/p} \big ) \ \ a.s.
$$
By  stationarity, this will hold true as soon as
\begin{equation} \label{conv1thgenep1}
\sum_{j \ge 1} \frac{\| \sup_{1 \leq k \leq 2^j} \big |  \sum_{\ell=1}^{k} \widetilde X_{j, \ell}\big |\|_2^2}{  2^{2j/p} j^{2-4/p}} < \infty \, .
\end{equation}
This is achieved by the following claim whose proof is given later.
\begin{claim} \label{lmaphi1} Assume that
\begin{equation}\label{cond1phi}
\sum_{k \geq 1} k^{-1/2} \phi^{1/2}_{1, {\bf Y}}(k) < \infty\, .
\end{equation}
Then \eqref{conv1thgenep1} holds.
\end{claim}
{\sc Proof of \eqref{conv2thgene}.} For any non negative integer $j$, let
\begin{eqnarray*} \label{0dec2}
\bar D_j := \sup_{ 1 \leq k \leq 2^{j}} \Big |  \sum_{\ell= 2^{j}+1}^{k+2^{j}}(\bar X_\ell - \bar d_{\ell}) \Big | =  \sup_{ 1 \leq k \leq 2^{j}} \Big |  \sum_{\ell= 2^{j}+1}^{k+2^{j}}(\bar X_{j,\ell} - \bar d_{j,\ell}) \Big |  \, .
\end{eqnarray*}
Following the beginning of the proof of (\ref{conv1thgene}), we get that $$\sup_{1\le k\le 2^N}| \bar S_k(f) -\bar M_k(f) |\le | \bar X_1-\bar d_1|+\sum_{j=0}^{N-1}\bar D_j \, .$$
Therefore  we infer that  \eqref{conv2thgene} will hold if
$$
\sup_{1\le k\le 2^{j}}
    \Big |\sum_{\ell= 2^{j}+1}^{k+2^{j}}(\bar X_{j,\ell} - \bar d_{j,\ell} ) \Big | = o \big (  j^{1 - 2/p}  \, 2^{j/p} \big )  \ \ a.s.
$$
By  stationarity, this will hold true as soon as
\begin{equation} \label{conv2thgenep1}
\sum_{j \ge 1}  \frac{\Big \|\sup_{1 \leq k \leq 2^j} \big |  \sum_{\ell=1}^{k} ( \bar X_{j,\ell} - \bar d_{j,\ell}) \big | \big \|_4^4}{  2^{4j/p} j^{4(1-2/p)}} < \infty \, .
\end{equation}
This is achieved by the following claim whose proof is given later.
\begin{claim}\label{lmaphi2}
Assume that
\begin{equation}\label{cond2phi}
\sum_{k \geq 1} k^{-1/8} \phi^{3/4}_{1, {\bf Y}}(k) < \infty\, .
\end{equation}
Then \eqref{conv2thgenep1} holds.
\end{claim}

\noindent {\sc Proof of \eqref{conv3thgene}.} Let $ a_n = n^{2/p} ( \log n)^{1-4/p}$. The result will follow from the following claim.
\begin{claim}\label{shlem} Assume that
\begin{equation} \label{condshao1}
\sum_{\ell \geq 1}  a_{\ell}^{-\nu} \E(|\bar d_{\ell}|^{2\nu}
) < \infty  \quad \mbox{ for some $1\le \nu \le 2$} \, ,
\end{equation}
and that
\begin{equation} \label{condshao2}
\sum_{i=1}^n\left({\mathbb E}_{i-1}(\bar d_i^{\, 2})-{\mathbb E}(\bar d_i^{\, 2})\right)
     =o ( a_n )\ \ a.s. 
\end{equation}
Then \eqref{conv3thgene} holds.
\end{claim}
It remains to show that \eqref{condshao1} and \eqref{condshao2} are satisfied.

\smallskip

\noindent {\sc Proof of \eqref{condshao1}}.  We show it  for $\nu =2$. Notice first that
$$
\sum_{\ell \geq 2}  a_{\ell}^{-2} \E(|\bar d_{\ell}|^{4}) = \sum_{j \geq 0}  \sum_{\ell= 2^j + 1}^{2^{j+1}}  a_{\ell}^{-2} \E(|\bar d_{j, \ell}|^{4}) \ll \sum_{j \geq 0}  \frac{1}{2^{4 j /p}  j^{2(1-4/p)}} \sum_{\ell= 2^j + 1}^{2^{j+1}} \E(|\bar d_{j, \ell}|^{4})  \, ,
$$
where $\bar d_{j , \ell}$ is defined in \eqref{defdj}. By stationarity and Lemma 5.1 in \cite{DDM},
\beq \label{p1lmaconvd4}
\Vert \bar d_{j, \ell} \Vert_4 \leq \sum_{\ell \geq 0} \Vert P_0 \left(\bar X_{j,\ell}\right) \Vert_4 \ll  \sum_{\ell \geq 0} (\ell+1)^{-1/4} \Vert \E_0 \left(\bar X_{j,\ell}\right) \Vert_4 \, .
\eeq
Using the arguments at the beginning of the proof of Claim \ref{convlqdjk}, we first observe that $(\bar X_{j,L, \ell})_{L \geq 1}$ defined in \eqref{defbartildeXjL} converges in ${\mathbb L}^4$ to $\bar X_{j,\ell}$. Hence
\beq \label{convevidentp4}
\Vert \E_0 \left(\bar X_{j,\ell}\right) \Vert_4 = \lim_{L \rightarrow \infty}\Vert \E_0 \left(\bar X_{j,L,\ell}\right) \Vert_4 \, .
\eeq
Next
$$
\Vert \E_0 \left(\bar X_{j,L,\ell}\right) \Vert_4 \leq  \sum_{k=1}^{\varphi(L)} |a_{k, \varphi(L)} | \Vert \E_0 \left( g_j \circ f_{k, \varphi(L)} (Y_{\ell}) \right )- \E \left( g_j \circ f_{k, \varphi(L)} (Y_{\ell}) \right ) \Vert_4 \, .
$$
Applying Lemma \ref{covaphi}, we then derive that
\[
\Vert \E_0 \left(\bar X_{j,L,\ell}\right) \Vert_4 \leq    2
  (2\phi_{1, {\bf Y}}(\ell ))^{3/4}  \sum_{k=1}^{\varphi(L)}
|a_{k, \varphi(L)}|  \Vert g_j \circ f_{k, \varphi(L)} (Y_0) \Vert_4 \, .\]
Therefore starting from \eqref{p1lmaconvd4}, we get that
\begin{align*}
\Vert \bar d_{j, \ell} \Vert_4  \ll \liminf_{L \rightarrow \infty} \sum_{k=1}^{\varphi(L)}
|a_{k, \varphi(L)}|  \Vert g_j \circ f_{k, \varphi(L)} (Y_0) \Vert_4  \sum_{\ell \geq 0} \ell^{-1/4} (\phi_{1, {\bf Y}}(k))^{3/4} \, .
\end{align*}
Since $\sum_{k=1}^{\varphi(L)}
|a_{k, \varphi(L)}| \leq 1$, Jensen's inequality leads to
$$
\Vert \bar d_{j, \ell} \Vert^4_4  \ll \liminf_{L \rightarrow \infty} \sum_{k=1}^{\varphi(L)}
|a_{k, \varphi(L)}|  \Vert g_j \circ f_{k, \varphi(L)} (Y_0)\Vert^4_4 \Big (  \sum_{\ell \geq 0} \ell^{-1/4} (\phi_{1, {\bf Y}}(\ell))^{3/4} \Big )^4\, .
$$
Hence, using condition \eqref{condDDM} and  Fubini, there exists a positive constant $C$ not depending on $L$ such that
\begin{align} \label{majnorme4d}
\sum_{j \geq 0}  & \frac{1}{2^{4 j /p}  j^{2(1-4/p)} } \sum_{\ell= 2^j + 1}^{2^{j+1}} \E(|\bar d_{j, \ell}|^{4}) \nonumber \\
&  \ll  \liminf_{L \rightarrow \infty} \sum_{k=1}^{\varphi(L)}
|a_{k, \varphi(L)}|  \sum_{j >0} \frac{2^j}{2^{4j/p} j^{2(p-4)/p}  }   \E \Big ( f^4_{k, \varphi(L)} (Y_0) {\bf 1}_{|f_{k, \varphi(L)} (Y_0)| \leq  2^{j/p} j^{-2/p}} \Big )   \nonumber \\
&\ll \liminf_{L \rightarrow \infty} \sum_{k=1}^{\varphi(L)}
|a_{k, \varphi(L)}|   \sum_{j >0}  \Big (  \frac{j^{2/p}  }{  2^{j/p} }  \Big )^{4-p}\E \Big ( f^4_{k, \varphi(L)} (Y_0) {\bf 1}_{|f_{k, \varphi(L)} (Y_0)| \leq  2^{j/p} j^{-2/p}} \Big )   \nonumber \\
& < C \liminf_{L \rightarrow \infty} \sum_{k=1}^{\varphi(L)}
|a_{k, \varphi(L)}|  \|f_{k, \varphi(L)} (Y_0) \|_p^p  \leq
C M^p\, .
\end{align}
This ends the proof of \eqref{condshao1}.

\medskip

\noindent {\sc Proof of (\ref{condshao2})}.  For any non negative integer $j$, let
\begin{eqnarray*}
A_j := \sup_{ 1 \leq \ell \leq 2^{j}} \big |  \sum_{i = 2^j +1}^{2^j + \ell} ({\mathbb E}_{i-1}({\bar d}_i{^{ 2}})-{\mathbb E}({\bar d}_i{^{ 2}}))\big | =  \sup_{ 1 \leq \ell \leq 2^{j}} \big |  \sum_{i = 2^j +1}^{2^j + \ell} ({\mathbb E}_{i-1}({\bar d}^{\, 2}_{j,i})-{\mathbb E}({\bar d}^{\, 2}_{j,i}) )\big |  \, ,
\end{eqnarray*}
where $\bar d_{j , i }$ is defined in \eqref{defdj}.
Let $N \in {\mathbb N}^*$ and let $k \in ]1, 2^N]$. We first notice that $A_j \geq |  \sum_{\ell= 2^{j}+1}^{2^{j+1}}({\mathbb E}_{i-1}({\bar d}^{\, 2}_{j,i})-{\mathbb E}({\bar d}^{\, 2}_{j,i}) )|$, so if $K$ is the  integer such that $2^{K-1} < k \leq  2^{K}$, then
$$\big  | \sum_{i=1}^k ({\mathbb E}_{i-1}({\bar d}^{\, 2}_{i})-{\mathbb E}({\bar d}^{\, 2}_{i}) ) \big | \leq |{\mathbb E}_{0}({\bar d}^{\, 2}_{1})-{\mathbb E}({\bar d}^{\, 2}_{1}) |+\sum_{j=0}^{N-1} A_j\, .$$
Consequently since $K \leq N$,
\begin{equation} \label{1dec} \sup_{1\le k\le 2^N}\big  | \sum_{i=1}^k ({\mathbb E}_{i-1}({\bar d}^{\, 2}_{i})-{\mathbb E}({\bar d}^{\, 2}_{i}) ) \big | \le |{\mathbb E}_{0}({\bar d}^{\, 2}_{1})-{\mathbb E}({\bar d}^{\, 2}_{1}) |+\sum_{j=0}^{N-1}A_j \, .
    \end{equation}
Therefore to prove \eqref{condshao2}, it is enough to show that
$$
\sup_{1\le k\le 2^{j}}
    \Big |\sum_{\ell= 2^{j}+1}^{k+2^{j}}({\mathbb E}_{i-1}({\bar d}^{\, 2}_{j,i})-{\mathbb E}({\bar d}^{\, 2}_{j,i}) )\Big | \leq  2^{2j/p} j^{1-4/p} \ \ a.s.
$$
In particular, it suffices to prove that
\begin{equation} \label{condshao2p2}
\sum_{j \ge 1} \frac{ \Big \Vert \sup_{1 \leq k \leq 2^j} \big |  \sum_{i=1}^{k} ({\mathbb E}_{i-1}({\bar d}^{\, 2}_{j,i})-{\mathbb E}({\bar d}^{\, 2}_{j,i}) ) \big |  \Big \Vert_2^2}{2^{4j/p} j^{2(1-4/p)} }  < \infty \, .
\end{equation}
Observe now that, by Claim \ref{convlqdjk} and condition \eqref{condDDM}, $(\bar d_{j,i})_{i \in {\mathbb Z}}$ forms a stationary triangular sequence of martingale differences in ${\mathbb L}^{4}$. Let
$$
\bar M_{j,n}:= \sum_{i=1}^n \bar d_{j,i} \, .
$$
Applying Proposition 2.3 in \cite{PU}  and using the martingale property of the sequence $(\bar M_{j,n})_{n \geq 1}$, we get that
$$
\Big \Vert \sup_{1 \leq k \leq 2^j} \big |  \sum_{i=1}^{k} ({\mathbb E}_{i-1}({\bar d}^{\, 2}_{j,i})-{\mathbb E}({\bar d}^{\, 2}_{j,i}) ) \big |  \Big \Vert_2^2\ll  2^j \| {\bar d}^{\, 2}_{j,1} \|^{2}_{2}  + 2^j \left ( \sum_{k=0}^{j-1} \frac{\|\E_0( \bar M_{j,2^k}^2) - \E(  \bar M_{j,2^k}^2) \|_{2}}{2^{k/2}}  \right )^{{2}} \, .
$$
Noticing that $\| d^2_{j,1} \|^2_{2} = \| d_{j,1} \|^4_{4}$, and using the computations made in \eqref{majnorme4d}, we get that
$$
\sum_{j \ge 1} \frac{ 2^j \| {\bar d}^{\, 2}_{j,1} \|^{2}_{2} }{2^{4j/p} j^{2(1-4/p)} } < \infty \, .
$$
Therefore to prove \eqref{condshao2p2}, it remains to prove that
$$
\sum_{j\ge 1} \frac{ 2^j  }{2^{4j/p} j^{2(1-4/p)} }\Big ( \sum_{k=0}^{j-1} 2^{-k/2} \|\E_0(\bar M_{j,2^k}^2) - \E(\bar M_{j,2^k}^2) \|_{2} \Big )^{2}<\infty \, .
$$
Setting $\bar S_{j,n} = \sum_{i=1}^n \bar X_{j,i}$ and $\bar R_{j,n}=\bar S_{j,n}-\bar M_{j,n}$, according to the proof of Theorem 2.3 in \cite{DDM}, this will hold provided that
\begin{equation} \label{condshao2p3}
\sum_{j\ge 1} \frac{ 2^j  }{2^{4j/p} j^{2(1-4/p)} }\Big ( \sum_{k=0}^{j-1} 2^{-k/2} \|\E_0(\bar S_{j,2^k}^2) - \E(\bar S_{j,2^k}^2) \|_{2} \Big )^{2}< \infty \, ,
\end{equation}
\begin{equation} \label{condshao2p4}
\sum_{j\ge 1} \frac{ 2^j  }{2^{4j/p} j^{2(1-4/p)} }\Big ( \sum_{k=0}^{j-1} 2^{-k/2} \| \bar R_{j,2^k} \|_{4}^{2} \Big )^2< \infty  \, ,
\end{equation}
and
\begin{equation} \label{condshao2p5}
\sum_{j\ge 1} \frac{ 2^j  }{2^{4j/p} j^{2(1-4/p)} }\Big ( \sum_{k=0}^{j-1} 2^{k/2}  \sum_{\ell \geq 2^k} \|P_0(\bar X_{j , \ell}) \|_{2}  \Big )^{2}< \infty\, .
\end{equation}
According to the proof of Claim \ref{lmaphi2},
$$
 \| \bar R_{j,2^k} \|_{4} \ll  \big ( \sum_{\ell\geq 1} (\phi_{1, {\bf Y}}(\ell))^{3/4}  \big ) \liminf_{L \rightarrow \infty}  \sum_{i=1}^{\varphi(L)}
|a_{i, \varphi(L)}|  \Vert g_j \circ f_{i, \varphi(L)} (Y_0) \Vert_4 \, .$$
\label{change} Therefore, since $\sum_{i=1}^{\varphi(L)}
|a_{i, \varphi(L)}| < 1$, using Jensen's inequality and condition \eqref{condDDM},
\begin{align*}
\sum_{j\ge 1} &\frac{ 2^j  }{2^{4j/p} j^{2(1-4/p)} }\Big ( \sum_{k=0}^{j-1} 2^{-k/2} \| \bar R_{j,2^k} \|_{4}^{2} \Big )^2 \\
& \ll \liminf_{L \rightarrow \infty}  \sum_{i=1}^{\varphi(L)}
|a_{i, \varphi(L)}|  \sum_{j\ge 1} \frac{ 2^j  }{2^{4j/p} j^{2(1-4/p)} } \E \Big ( f^4_{i, \varphi(L)} (Y_0) {\bf 1}_{|f_{i, \varphi(L)} (Y_0)| \leq  2^{j/p} j^{-2/p}} \Big ) \, ,
\end{align*}
which together with an application of Fubini show \eqref{condshao2p4}. We turn to the proof of \eqref{condshao2p5}. According to Lemma 5.1 in \cite{DDM},
$$
2^{k/2} \sum_{\ell \geq 2^k} \|P_0(\bar X_{j , \ell}) \|_{2} \ll \sum_{\ell \geq 2^{k-1}}  \|\E_0(\bar X_{j , \ell}) \|_{2} \, .
$$
Since $\|\E_0(\bar X_{j , \ell}) \|_{2} \leq \liminf_{L \rightarrow \infty}  \sum_{i=1}^{\varphi(L)}
|a_{i, \varphi(L)}|  \|\E_0(g_j \circ f_{i, \varphi(L)} (Y_\ell)) - \E(g_j \circ f_{i, \varphi(L)} (Y_\ell)) \|_{p}$, by Lemma \ref{covaphi} and using the fact that $ \sum_{i=1}^{\varphi(L)}
|a_{i, \varphi(L)}|  \leq 1$, we derive that
$$\|\E_0(\bar X_{j , \ell}) \|_{2} \ll M ( \phi_{1, {\bf Y}}(\ell))^{(p-1)/p}\, .$$
Since $p <4$, we then infer that \eqref{condshao2p5} will hold if $\sum_{\ell \geq 1} \phi_{1, {\bf Y}}(\ell))^{(p-1)/p} < \infty $ which is satisfied if \eqref{condDDM} is.

It remains to show that \eqref{condshao2p3} holds true. Since
$$
 \Vert \E_0(\bar S_{j,2^k}^2) - \E(\bar S_{j,2^k}^2) \Vert_{2}
 \leq 2 \sum_{m=1}^{2^k}\sum_{\ell=0}^{2^k-m}\Vert \E_0 (\bar X_{j,m} \bar X_{j,m+\ell}) - \E (\bar X_{j,m} \bar X_{j,m+\ell}) \Vert_{2} \, , $$ we infer that
there exists $C>0$ such that
\begin{align*}\label{b1}
\sum_{k=0}^{j-1} & 2^{-k/2} \|\E_0(\bar S_{j,2^k}^2) - \E(\bar S_{j,2^k}^2) \|_{2}  \nonumber \\
& \leq C \sum_{m=1}^{2^j}\sum_{\ell=0}^{2^j}
\frac{1}{ (m+\ell)^{1/2}} \Vert \E_0 (\bar X_{j,m} \bar X_{j,m+\ell}) - \E (\bar X_{j,m} \bar X_{j,m+\ell}) \Vert_{2}\nonumber \\
& \leq C \sum_{m=1}^{2^j}\sum_{\ell=0}^{[m^{\gamma}]}
\frac{1}{ (m+\ell)^{1/2}}  \gamma_j({\mathcal F}_0, m,\ell)  +C \sum_{m=1}^{2^j} \sum_{\ell =[m^{\gamma}]}^{2^j}
\frac{1}{ (m+\ell)^{1/2}}  \gamma_j({\mathcal F}_0, m,\ell)  \, .
\end{align*}
where
$$
\gamma_j({\mathcal F}_0, m,\ell) :=\Vert \E_0 (\bar X_{j,m} \bar X_{j,m+\ell}) - \E (\bar X_{j,m} \bar X_{j,m+\ell}) \Vert_{2} \, \text{ and } \, \gamma \in (0,1]\, .
$$
Bounding up $\gamma_j({\mathcal F}_0, m,\ell) $ in two ways as done in the proof of Corollary 2.1 in \cite{DDM}, we infer that, for any $\gamma\in (0,1]$ (to be chosen later), there exists
a positive constant $B$ such that
\begin{equation}\label{borneDDM}
\Big ( \sum_{k=0}^{j-1} 2^{-k/2} \|\E_0(\bar S_{j,2^k}^2) - \E(\bar S_{j,2^k}^2) \|_{2} \Big )^{2}\leq B I_1^2 +B I_2^2
\end{equation}
where
\begin{align*}
I_1 &= \sum_{m=1}^{2^j} \frac{m^{\gamma}}{m^{1/2}} \sup_{\ell  \geq k \geq m}
\Vert \E_0 (\bar X_{j,\ell} \bar X_{j,k}) - \E (\bar X_{j,\ell} \bar X_{j,k}) \Vert_{2}  \\
I_2 &= \Big(\sum_{k=1}^{2^j} \frac{k^{1/(2\gamma)}}{k^{1/4}}
\|\E_0( \bar X_{j,k})\|_4\Big)^2\, .
\end{align*}
We shall proceed by using some arguments developed in \cite{DeGoMe} to get their bound (5.7). Arguing as in the proof of Claim \ref{lmaphi2}, we obtain that
\begin{equation} \label{I2}
I_2 \leq 8 \sqrt 2 \Big(  \sum_{k>0} \frac{k^{1/(2\gamma)}}{k^{1/4}}
\phi_{1, {\bf Y}}(k)^{3/4} \Big)^2 \liminf_{L \rightarrow \infty}
\Big(\sum_{\ell =1}^{\varphi(L)}
|a_{\ell, \varphi(L)}|
\E(f_{\ell, \varphi(L)}^4 (Y_0) {\bf 1}_{|f_{\ell, \varphi(L)} (Y_0)| \leq  c(j)})\Big)^{1/2} \, .
\end{equation}
We bound now $I_1$. According to the proof of Claim \ref{convlqdjk}, $\bar X_{j,\ell}$ is the limit in $\LL^4$
of $ (\bar X_{j,L, \ell})_L$, where $\bar X_{j,L, \ell}$ is defined in
 \eqref{defbartildeXjL}. Therefore,
$$
\Vert \E_0 (\bar X_{j,\ell} \bar X_{j,k}) - \E (\bar X_{j,\ell} \bar X_{j,k}) \Vert_{2} =  \lim_{L \rightarrow \infty}\Vert \E_0 (\bar X_{j,L, \ell} \bar X_{j,L, k}) - \E (\bar X_{j,L,\ell} \bar X_{j,L, k}) \Vert_{2}\, .
$$
Applying Lemma \ref{covaphi}, for $\ell  \geq k \geq m$,
\begin{multline*}
\Vert \E_0 (\bar X_{j,L, \ell} \bar X_{j,L, k}) - \E (\bar X_{j,L,\ell} \bar X_{j,L, k}) \Vert_{2}\\
 \leq 16  \phi_{2, {\bf Y}}(m)^{1/2} \sum_{\ell =1}^{\varphi(L)}
|a_{\ell, \varphi(L)}|  \|g_j \circ f_{\ell, \varphi(L)}(Y_0)\|_4 \sum_{k =1}^{\varphi(L)}
|a_{k, \varphi(L)}|
 \|g_j \circ f_{k, \varphi(L)} (Y_0)\|_4 \, .
\end{multline*}
It follows that
\begin{equation}\label{I1}
 I_1 \leq \Big(16 \sum_{m=1}^{2^j} \frac{m^{\gamma}}{m^{1/2}} \phi_{2, {\bf Y}}(m)^{1/2}\Big)
\liminf_{L \rightarrow \infty}
\Big(\sum_{\ell =1}^{\varphi(L)}
|a_{\ell, \varphi(L)}|
\E(f_{\ell, \varphi(L)}^4 (Y_0) {\bf 1}_{|f_{\ell, \varphi(L)} (Y_0)| \leq  c(j)})\Big)^{1/2}  \, .
\end{equation}
Let $\gamma=1/\sqrt 3$. If the condition \eqref{condDDM} holds, then
\[
\sum_{k>0} \frac{k^{{\sqrt 3}/2}}{k^{1/4}}
\phi_{1, {\bf Y}}(k)^{3/4}< \infty \quad \text{and} \quad
\sum_{m>0} \frac{m^{1/\sqrt 3}}{m^{1/2}} \phi_{2, {\bf Y}}(m)^{1/2} < \infty.
\]

We infer  from \eqref{borneDDM}, \eqref{I2} and \eqref{I1} that, if
\eqref{condDDM} holds, there exists a positive constant $C_4(\phi)$ such
that
\begin{equation}\label{borne}
\Big ( \sum_{k=0}^{j-1} 2^{-k/2} \|\E_0(\bar S_{j,2^k}^2) - \E(\bar S_{j,2^k}^2) \|_{2} \Big )^{2} \leq C_4(\phi)
\liminf_{L \rightarrow \infty}
\sum_{\ell =1}^{\varphi(L)}
|a_{\ell, \varphi(L)}|
\E(f_{\ell, \varphi(L)}^4 (Y_0) {\bf 1}_{|f_{\ell, \varphi(L)} (Y_0)| \leq  c(j)}) \, .
\end{equation}
Using this last bound, Fatou's lemma together with Fubini, we then infer that \eqref{condshao2p3} holds true. This ends the proof of \eqref{conv3thgene}.
$\square$

\medskip
\noindent{\bf End of the proof Theorem \ref{asipgenephi} for $p \in ]2,4[$.}
To finish the proof of Theorem \ref{asipgenephi} for $p \in ]2,4[$, it remains to prove our claims.

\medskip

\noindent{\bf Proof of Claim \ref{convlqdjk}.} Notice that $(\bar X_{j,L, \ell})$ converges in $\LL^p$ to $\bar X_{j,\ell}$ and $\Vert \bar X_{j,L, \ell} \Vert_{\infty} \leq 2 c(j)$.
Therefore, we infer that $\Vert \bar X_{j, \ell} \Vert_{\infty} \leq 2 c(j)$. It follows that $\bar X_{j, \ell}  =  \lim_{L \rightarrow \infty} \bar X_{j,L,\ell}$ in $\LL^q$ for any $q \in [1, \infty [$. Next, by Lemma \ref{covaphi} and the fact that $\sum_{k=1}^{\varphi(L)}|a_{k, \varphi(L)}| \leq 1$, we get that
\begin{align*}
\Vert \E_0 \left(\bar X_{j,L,\ell}\right) \Vert_q & \leq \sum_{k=1}^{\varphi(L)} |a_{k, \varphi(L)} | \Vert \E_0 \left( g_j \circ f_{k, \varphi(L)} (Y_{\ell}) \right )- \E \left( g_j \circ f_{k, \varphi(L)} (Y_{\ell}) \right ) \Vert_{\infty} \\
& \leq 4
  \phi_{1, {\bf Y}}(\ell)  \sum_{k=1}^{\varphi(L)}
|a_{k, \varphi(L)}|  \Vert g_j \circ f_{k, \varphi(L)} (Y_0)\Vert_{\infty} \leq 4 c(j)
  \phi_{1, {\bf Y}}(\ell)  \, .
\end{align*}
Therefore $\Vert \E_0 \left(\bar X_{j,\ell}\right) \Vert_q \leq 4 c(j)
  \phi_{1, {\bf Y}}(\ell)$. This proves the claim. $\square$

\medskip

\noindent{\bf Proof of Claim \ref{lmaphi1}.} Let $\widetilde S_{j,L,k} = \sum_{\ell=1}^{k} \widetilde  X_{j, L, \ell}$ where $\widetilde  X_{j, L, \ell}$ is defined by \eqref{defbartildeXjL}. Recall that $\widetilde X_{j,\ell}$ is the limit in ${\mathbb L}^p$ of $\widetilde X_{j,L,\ell}$. Clearly, to prove \eqref{conv1thgenep1}, it suffices to prove that, there exists
some positive constant $K$, such that for any positive integer $L$,
\begin{equation} \label{but2bis}
\sum_{j \ge 1}\frac{\| \sup_{1 \leq k \leq 2^j} \big |  \widetilde S_{j,L,k}
\big |\|_2^2}{  2^{2j/p} j^{2-4/p}} < K \, .
\end{equation}
To prove \eqref{but2bis}, we use the maximal inequality of Peligrad and Utev (2006). Therefore, by stationarity,
\begin{align} \label{ineMS}
\big \Vert \max_{1 \leq k \leq 2^j} \big | \widetilde  S_{j,L,k} \big | \big \Vert^2_2 & \ll 2^j \Vert \widetilde  X_{j, L, 0} \Vert^2_2 + 2^j \Big ( \sum_{\ell=0}^{j} 2^{-\ell/2} \Vert \E_0 ( \widetilde  S_{j,L,2^{\ell}} )\Vert_2 \Big )^2 \nonumber \\
& \ll 2^j \Vert \widetilde X_{j, L, 0} \Vert^2_2 +  2^j \Big ( \sum_{k=1}^{2^j} k^{-1/2} \Vert \E_0 ( \widetilde  X_{j,L,k} )\Vert_2 \Big )^2 \, .
\end{align}
Notice  that $
\Vert \widetilde X_{j, L, 0} \Vert^2_2 \leq 4 \Big ( \sum_{\ell=1}^{\varphi(L)} a_{\ell, \varphi(L)}
\Vert \widetilde g_j\circ  f_{ \ell,\varphi( L)} (Y_0) \Vert_2 \Big )^2
$. Therefore, since $\sum_{\ell=1}^L
|a_{\ell, L}| \leq 1$,  by Jensen's inequality,
\beq \label{1ms}
\Vert \widetilde X_{j, L, 0} \Vert^2_2 \leq 4  \sum_{\ell=1}^{\varphi(L)} | a_{\ell, \varphi(L)} |
\Vert  \widetilde g_j\circ  f_{ \ell,\varphi( L)} (Y_0)  \Vert_2^2 \, .
\eeq
Now
\[
\Vert \E_0 ( \widetilde  X_{j,L,k} )\Vert_2
\leq  \sum_{\ell=1}^{\varphi(L)} | a_{\ell, \varphi(L)} |   \|\E_0( \widetilde g_j\circ  f_{ \ell,\varphi( L)} (Y_k) )-\E( \widetilde g_j\circ  f_{ \ell,\varphi( L)} (Y_k) ) \|_2 \, .
\]
Applying Lemma \ref{covaphi},
$ \|\E_0( \widetilde g_j\circ  f_{ \ell,\varphi( L)} (Y_k) )-\E( \widetilde g_j\circ  f_{ \ell,\varphi( L)} (Y_k) ) \|_2 \leq 2(2\phi_{1, {\bf Y}}(k))^{1/2}\| \widetilde g_j\circ  f_{ \ell,\varphi( L)} (Y_0)\|_2$. Hence by Jensen's inequality,
\begin{equation}\label{2ms}
\Big ( \sum_{k=1}^{2^j} k^{-1/2} \Vert \E_0 ( \widetilde  X_{j,L,k} )\Vert_2 \Big )^2  \leq 8  \Big ( \sum_{k\geq 1} k^{-1/2}\phi^{1/2}_{1, {\bf Y}}(k) \Big )^2  \sum_{\ell=1}^{\varphi(L)} | a_{\ell, \varphi(L)} |  \| \widetilde g_j\circ  f_{ \ell,\varphi( L)} (Y_0)\|_2^2 \, .
\end{equation}
Therefore, using \eqref{ineMS} together with the upper bounds \eqref{1ms} and \eqref{2ms}, we derive that
\begin{align*}
\sum_{j \ge 1}\frac{\| \max_{1 \leq k \leq 2^j} \big |  \widetilde S_{j,L,k}
\big |\|_2^2}{  2^{2j/p} j^{2-4/p}}
& \ll  \frac{2^{j(p-2)/p}}{j^{2(p-2)/p} } \Big (  \sum_{k\geq 1} k^{-1/2}\phi^{1/2}_{1, {\bf Y}}(k) \Big )^2\sum_{\ell=1}^{\varphi(L)} | a_{\ell, \varphi(L)} |  \| \widetilde g_j\circ  f_{ \ell,\varphi( L)} (Y_0)\|_2^2 \, .
\end{align*}
Now, via Fubini, there exists a positive constant $C$ not depending on $L$ such that
\begin{align*}
\sum_{j >0} \frac{2^{j(p-2)/p}}{j^{2(p-2)/p} }   \| \widetilde g_j\circ  f_{ \ell,\varphi( L)}(Y_0)\|_2^2 & = \sum_{j >0} \frac{2^{j(p-2)/p}}{j^{2(p-2)/p} }  \E \Big ( f^2_{\ell, \varphi(L)} (Y_0) {\bf 1}_{|f_{\ell, \varphi(L)} (Y_0)| > 2^{j/p} j^{-2/p}} \Big )  \\
& < C \|f_{\ell, \varphi(L)} (Y_0) \|_p^p  \leq
C M^p\, .
\end{align*}
Using condition \eqref{cond1phi} and the fact that $\sum_{\ell=1}^{\varphi(L)} | a_{\ell, \varphi(L)} | <1$, \eqref{but2bis} follows. This ends the proof of the
claim. $\square$

\medskip

\noindent{\bf Proof of Claim \ref{lmaphi2}.} For any positive integer $k$, let $ \bar S_{j,k} =\sum_{\ell=1}^{k} \bar X_{j, \ell}$, $ \bar M_{j,k} = \sum_{\ell=1}^{k} \bar d_{j, \ell}$ and
\beq \label{defRbarjL}\bar R_{j,k}  = \bar S_{j,k} -  \bar M_{j,k} \, .\eeq
By stationarity, and since for any $k \geq 1 $, $\E_0(  \bar M_{j,k} )= 0 $, according to Corollary 3 in Merlev\`ede and Peligrad \cite{MP},
\begin{align*}
\Vert \sup_{1 \leq k \leq 2^j} \text{ }| \bar R_{j,k}|\text{ }\Vert_{4} & \leq \Vert
\bar R_{j,2^j}\Vert_{4}+2^{j/4}\sum_{l=0}^{j-1}2^{-l/4} \Vert \E_0(  \bar S_{j,2^l} )\Vert_{4} \nonumber  \\
& \ll \Vert
\bar R_{j,L,2^j}\Vert_{4}+2^{j/4}\sum_{\ell=1}^{2^j} \ell^{-1/4} \Vert \E_0(  \bar X_{j,\ell} )\Vert_{4} \, .
\end{align*}
Now according to item 2 of Proposition 2.1 in \cite{DDM} applied with $N= 2^{2j+1}$, we get that
$$
\Vert
\bar R_{j,2^j}\Vert_{4} \ll \sum_{\ell=1}^{2^{2j}} \Vert \E_0(  \bar X_{j,\ell} )\Vert_{4} + 2^{j/2} \sum_{\ell \geq 2^{2j}  } \Vert P_0(  \bar X_{j,\ell} )\Vert_{4} \, .
$$
Using then Lemma 5.1 in \cite{DDM}, it follows that$$
\Vert
\bar R_{j,2^j}\Vert_{4} \ll \sum_{\ell=1}^{2^{2j +1}} \Vert \E_0(  \bar X_{j,\ell} )\Vert_{4} + 2^{j/2} \sum_{\ell \geq 2^{2j} } k^{-1/4}\Vert \E_0(  \bar X_{j,\ell} )\Vert_{4} \, .
$$
So overall,
\begin{align} \label{p1lma2}
\Vert \sup_{1 \leq k \leq 2^j} & \text{ }| \bar R_{j,k}|\text{ }\Vert_{4}  \nonumber \\
& \leq \sum_{\ell=1}^{2^{2j +1}} \Vert \E_0(  \bar X_{j,\ell} )\Vert_{4}  + 2^{j/2} \sum_{\ell \geq 2^{2j}  } k^{-1/4}\Vert \E_0(  \bar X_{j,\ell} )\Vert_{4} +2^{j/4}\sum_{\ell=1}^{2^j} \ell^{-1/4} \Vert \E_0(  \bar X_{j,\ell} )\Vert_{4} \nonumber \\
& \leq 2^{j/4}\sum_{\ell=1}^{2^j} \ell^{-1/8} \Vert \E_0(  \bar X_{j,\ell} )\Vert_{4}  + 2^{j/2} \sum_{\ell \geq 2^{2j}  } \ell^{-1/4}\Vert \E_0(  \bar X_{j,\ell} )\Vert_{4} \, .
\end{align}
We handle now the quantity $\Vert \E_0(  \bar X_{j,\ell} )\Vert_{4} $. We first observe that by Lemma \ref{covaphi},
$$
\Vert \E_0 \left(\bar X_{j,L,\ell}\right) \Vert_4 \leq 4  (\phi_{1, {\bf Y}}(\ell))^{3/4}\sum_{k=1}^{\varphi(L)}
|a_{k, \varphi(L)}|  \Vert g_j \circ f_{k, \varphi(L)} (Y_0) \Vert_4 \, ,$$
where $X_{j,L,\ell}$ is defined in \eqref{defbartildeXjL}. Taking into account \eqref{convevidentp4}, \eqref{p1lma2} and the condition \eqref{cond2phi}, we then infer that
\[
\Vert \sup_{1 \leq k \leq 2^j}  \text{ }| \bar R_{j,k}|\text{ }\Vert_{4}^4  \ll
\frac{2^j}{ 2^{4j/p} j^{4(p-2)/p}  } \liminf_{L \rightarrow \infty} \Big (  \sum_{k=1}^{\varphi(L)}
|a_{k, \varphi(L)}|  \Vert g_j \circ f_{k, \varphi(L)} (Y_0) \Vert_4 \Big )^{4} \, .
\]
Since, $\sum_{\ell=1}^{\varphi(L)} | a_{\ell, \varphi(L)} | <1$, by Jensen's inequality
$$
\Vert \sup_{1 \leq k \leq 2^j}  \text{ }| \bar R_{j,k}|\text{ }\Vert_{4}^4
\ll \frac{2^j}{ 2^{4j/p} j^{4(p-2)/p}  } \liminf_{L \rightarrow \infty}  \sum_{k=1}^{\varphi(L)}
|a_{k, \varphi(L)}|  \Vert g_j \circ f_{k, \varphi(L)} (Y_0) \Vert_4^{4} \, .
$$
Now, via Fubini, there exists a positive constant $C$ not depending on $L$ such that
\begin{align*}
\sum_{j >0} & \frac{2^j}{2^{4j/p} j^{4(p-2)/p}  }   \Vert g_j \circ f_{k, \varphi(L)} (Y_0) \Vert_4^{4}  = \sum_{j >0} \frac{2^j}{ 2^{4j/p} j^{4(p-2)/p}  }   \E \Big ( f^4_{\ell, \varphi(L)} (Y_0) {\bf 1}_{|f_{\ell, \varphi(L)} (Y_0)| \leq  2^{j/p} j^{-2/p}} \Big )  \\
&\leq \ \sum_{j >0}  \Big (  \frac{2^{j/p} }{ j^{2/p}  }  \Big )^{p-4}\E \Big ( f^4_{\ell, \varphi(L)} (Y_0) {\bf 1}_{|f_{\ell, \varphi(L)} (Y_0)| \leq  2^{j/p} j^{-2/p}} \Big )   < C \|f_{\ell, \varphi(L)} (Y_0) \|_p^p  \leq
C M^p\, ,
\end{align*}
which combined with the fact that  $\sum_{\ell=1}^{\varphi(L)} | a_{\ell, \varphi(L)} | <1$ prove the claim.  $\square$

\medskip

\noindent {\bf Proof of Claim \ref{shlem}.} We first prove that
\beq \label{limbarMn}
\sigma^2 = \lim_{n \rightarrow \infty}n^{-1} \Vert \bar M_n(f) \Vert_2^2 \, .
\eeq
Recall that \eqref{condDDM} entails in particular that   $\sum_{k \geq 0} \Vert P_0 (X_k) \Vert_2 < \infty$. We then define $d_0 = \sum_{i \geq 0} P_0 (f (Y_0) \circ \theta^i)$ and for any integer $\ell$, $d_{\ell}= d_0 \circ \theta^{\ell}$. Let $M_n (f)= \sum_{\ell=1}^n d_{\ell} $. Since $\sum_{k \geq 0} \Vert P_0 (X_k) \Vert_2 < \infty$, using item 2 of Theorem 1 in \cite{W},
\begin{equation} \label{app2**}
\Vert S_n (f) - M_n (f) \Vert_2 = o (\sqrt{n}) \, .
\end{equation}
Since $\sigma^2 = \lim_{n \rightarrow \infty} n^{-1} \Vert S_n(f) \Vert_2^2$, it follows from \eqref{app2**} and stationarity that $\sigma^2 =n^{-1} \Vert M_n(f) \Vert_2^2= \E (d_0^2)$. We show now that
\begin{equation} \label{condshao3}
\Vert \bar M_n (f)- M_n (f) \Vert_2  =o(n^{1/2}) \, .
\end{equation}
Let $N$ be the positive integer such that $ 2^{N-1} < n  \leq 2^N$. By the martingale property of  $\bar M_n (f)- M_n (f)$ and stationarity, we have that
\begin{equation} \label{estnorm2}
\Vert \bar M_n (f)- M_n (f)  \Vert^2_2 = \sum_{\ell = 1}^n \E ( ( \bar d_{\ell} - d_{\ell})^2) \leq \E ( (\bar d_{1} - d_{1})^2) +  \sum_{j=0}^{N-1}2^j \E ( (d_{ 1}-  \bar d_{j,1}  )^2) \, .
\end{equation}
But $ d_{ 1}-  \bar d_{j,1} = \sum_{\ell \geq 1} P_1  (  \widetilde X_{j,\ell}  ) $. Then, by Lemma 5.1 in \cite{DDM} (see also the proof of Corollary 2 in \cite{PU06}),
$$
\Vert d_{ 1}-  \bar d_{j,1} \Vert_2 \leq \sum_{\ell \geq 0} \Vert P_0  (  \widetilde X_{j,\ell}  ) \Vert_2 \ll \sum_{\ell \geq 0} (\ell+1)^{-1/2} \Vert \E_0  (  \widetilde X_{j,\ell}  ) \Vert_2 \, .
$$
Let $\widetilde X_{j,L, \ell}$ be defined in \eqref{defbartildeXjL}. Applying Lemma \ref{covaphi},
\[
\Vert \E_0 (\widetilde X_{j,L, \ell}) \Vert_2 \leq    2 \sqrt{2}
 ( \phi_{1, {\bf Y}}(\ell))^{1/2}  \sum_{k=1}^{\varphi(L)}
|a_{k, \varphi(L)}|  \Vert \widetilde g_j \circ f_{k, \varphi(L)} (Y_0) \Vert_{2}  \, .
\]
Since $\sum_{k=1}^{\varphi(L)}
|a_{k, \varphi(L)}| \leq 1$ and $f_{k, \varphi(L)} $ belongs to  $\Mon_p(M,P_{Y_0})$, it follows that
\[
\Vert \E_0 (\widetilde X_{j,L,\ell}) \Vert_2 \leq    2 \sqrt{2} M^{p/2} (c(j))^{1-p/2}
  (\phi_{1, {\bf Y}}(k))^{1/2}\, .
\]
Since $\Vert \E_0  (  \widetilde X_{j,\ell}  ) \Vert_2 = \lim_{L \rightarrow \infty} \Vert \E_0 (\widetilde X_{j,L,\ell}) \Vert_2$  and, by condition \eqref{condDDM},  $\sum_{k \geq 1} k^{-1/2} \phi^{1/2}_{1, {\bf Y}}(k) < \infty$, we get overall that
$$
\Vert \bar M_n (f)- M_n (f)  \Vert^2_2  \ll   \sum_{j=0}^{N-1}2^j (c(j))^{2-p} \ll 2^{2N/p} N^{2(1-2/p)} \, ,
$$
proving \eqref{condshao3}. Combining the fact that $\sigma^2 =n^{-1} \Vert M_n(f) \Vert_2^2$ with \eqref{condshao3}, it follows that \eqref{limbarMn} holds.
\medskip

Then, according to Theorem 2.1 in \cite{shao93},
we see that, enlarging our probability space  if
necessary, one may find a sequence $(\bar Z_\ell)_{\ell
\geq 1}$ of independent  Gaussian random variables with zero mean and
variance $\E( \bar Z_\ell)^2= \E ( \bar d_{\ell})^2= ( \bar \sigma_{\ell} )^2$ such that
\begin{equation} \label{asipmart1}
\sup_{1\leq k \leq n} \Big| \bar M_k - \sum_{\ell=1}^k \bar Z_{\ell}\Big|  = o \big ( a_n^{1/2} ( \log n)^{1/2}  \big )
\text{ almost surely, as $n\rightarrow \infty$}.
\end{equation}
Let $(\delta_k)_{k \geq
1}$ be a sequence of iid\ Gaussian random variables with mean zero
and variance $\sigma^2$, independent of the sequence $(\bar Z_{\ell})_{\ell\geq
1}$. We now construct a sequence $(Z_{\ell})_{\ell \geq 1}$ as follows. If $\bar \sigma_{\ell}=0$, then $Z_\ell=\delta_\ell$, else $Z_\ell=
(\sigma/ \bar \sigma_{\ell}) \bar Z_{\ell}$. By construction, the $Z_\ell$'s are iid\ Gaussian random variables with
mean zero and variance $\sigma^2$. Let $G_\ell=Z_\ell-\bar Z_{\ell}$ and note that $(G_\ell)_{\ell \geq 1}$ is a
sequence of independent Gaussian random variables with mean zero and
variances $\Var(G_\ell)=(\sigma-\bar \sigma_{\ell})^2$. Assume that we can prove that
\beq \label{kolmogorovcond}
\sum_{n \geq 3} \frac{\E (G_n^2)}{a_n \log n} < \infty \, .
\eeq
Then by the Kolmogorov theorem (or Lemma \ref{reverse}), it will follow that the series $\sum_{n \geq 3} \frac{G_n}{(a_n \log n)^{1/2}}$ converges $\as$ Hence, Kronecker lemma will imply that $\sum_{\ell =1}^n G_{\ell} = o ((a_n \log n)^{1/2} )$ $\as$ Therefore starting from \eqref{asipmart1}, we will conclude that if \eqref{condshao1} and \eqref{condshao2} hold then \eqref{conv3thgene} does. Let us prove \eqref{kolmogorovcond}. With this aim, we first notice that
$$
\E (G_n^2) = \big (  \Vert d_n \Vert_2 -\Vert \bar d_n\Vert_2 \big )^2 \leq  \Vert d_n -\bar d_n\Vert^2_2 \, .
$$
Next
$$
\sum_{n \geq 3} \frac{\E (G_n^2)}{a_n \log n}  \leq \sum_{j \geq 1} \frac{1}{  2^{2j/p} j^{2-4/p}} \sum_{\ell= 2^{j}+1}^{2^{j+1}} \E ( (d_{ \ell}-  \bar d_{j,\ell}  )^2) =\sum_{j \geq 1} \frac{2^j}{  2^{2j/p} j^{2-4/p}}  \E ( (d_{ 1}-  \bar d_{j,1}  )^2) \, .
$$
Using the computations as done to prove \eqref{condshao3}, we infer that under condition \eqref{condDDM},
$$
\sum_{n \geq 3} \frac{\E (G_n^2)}{a_n \log n}  \ll \liminf_{L \rightarrow \infty}\sum_{j \geq 1} \frac{2^j}{  2^{2j/p} j^{2-4/p}}  \Big ( \sum_{k=1}^{\varphi(L)}
|a_{k, \varphi(L)}|  \Vert \widetilde g_j \circ f_{k, \varphi(L)} (Y_0)\Vert_2 \Big )^2  \, .
$$
Since $\sum_{k=1}^{\varphi(L)}
|a_{k, \varphi(L)}| \leq 1$, Jensen's inequality leads to
$$
\sum_{n \geq 3} \frac{\E (G_n^2)}{a_n \log n}  \ll \liminf_{L \rightarrow \infty}\sum_{j \geq 1} \frac{2^j}{  2^{2j/p} j^{2-4/p}} \sum_{k=1}^{\varphi(L)}
|a_{k, \varphi(L)}|   \Vert \widetilde g_j \circ f_{k, \varphi(L)} (Y_0)\Vert_2^2 \, .
$$
Hence, by Fubini theorem, there exists a positive constant $C$ not depending on $L$ such that
\begin{align*}
\sum_{n \geq 3}  \frac{\E (G_n^2)}{a_n \log n} & \ll \liminf_{L \rightarrow \infty}\sum_{j \geq 1} \frac{2^j}{  2^{2j/p} j^{2-4/p}} \sum_{k=1}^{\varphi(L)}
|a_{k, \varphi(L)}|   \Vert \widetilde g_j \circ f_{k, \varphi(L)} (Y_0)\Vert_2^2   \\
&  \ll  \liminf_{L \rightarrow \infty} \sum_{k=1}^{\varphi(L)}
|a_{k, \varphi(L)}|  \sum_{j \geq 1} \Big ( \frac{2^{j/p}}{  j^{2/p}} \Big )^{p-2}  \E \Big ( f^2_{k, \varphi(L)} (Y_0) {\bf 1}_{|f_{k, \varphi(L)} (Y_0)| > 2^{j/p} j^{-2/p}} \Big )    \\
& < C \liminf_{L \rightarrow \infty} \sum_{k=1}^{\varphi(L)}
|a_{k, \varphi(L)}|  \|f_{k, \varphi(L)} (Y_0) \|_p^p  \leq
C M^p\, .
\end{align*}
This ends the proof of \eqref{kolmogorovcond} and of Claim \ref{shlem}. $\square$



\subsection{Proof of Theorem \ref{asipgenephi} for $p=4$}
In this case, no truncation is needed. The beginning of the proof of Theorem \ref{asipgenephi} for $2<p<4$ also works for $p=4$. In
particular, if
$f\in \Monm_4(M,P_{Y_0})$, condition \eqref{condDDM} implies that $\sum_{\ell \geq 0} \Vert \E_0(X_{\ell} )\Vert_4< \infty$. Therefore the series $\sigma^2= \sum_{k \in {\mathbb Z}} {\rm Cov} ( X_0,X_k)$  converges absolutely and   $ \lim_{n \rightarrow \infty }n^{-1} \E (S^2_n(f)) = \sigma^2$. In addition the series $\sum_{k \geq 0}P_0\left( X_{k}\right)$ converges in ${\mathbb L}^{4}$ and the sequence  $(d_{\ell})_{\ell \in {\mathbb Z}}$ defined by:
\beq \label{defdlp=4}
d_{\ell} = \sum_{k \geq \ell}P_\ell\left( X_{k}\right)
\eeq
forms a  stationary sequence of martingale differences in ${\mathbb L}^4$ with respect to the non-decreasing sequence of  $\sigma$-algebras $({\mathcal F}_\ell)_{\ell \in {\mathbb Z}}$. Hence, setting for  every positive integer $n$,
$$  M_n(f):=\sum_{\ell=1}^n  d_\ell\ \ \mbox{and}\ \  S_n(f):=\sum_{\ell=1}^n  X_{\ell} \, ,$$the conclusion of Theorem \ref{asipgenephi} when $p=4$ will follow if we can prove that
\beq \label{conv1thgenep=4}
\sup_{1 \leq k \leq n } \big | S_n(f)-  M_n(f)\big | = O ( n^{1/4} (\log n)^{1/2} ( \log \log n)^{1/4}) \text{ almost surely,}
\eeq
and if, enlarging our probability space if necessary, there exists  a
    sequence $(Z_i)_{i \geq 0}$ of iid\ Gaussian random
    variables with mean zero and variance $\sigma^2$
such that
\begin{equation}\label{conv3thgenep=4}
\sup_{1 \leq k \leq n }\Big|\sum_{i=1}^{k} ( d_i
- Z_i) \Big| =O ( n^{1/4} (\log n)^{1/2} ( \log \log n)^{1/4})\,
  \text{almost surely.}
\end{equation}
To prove \eqref{conv1thgenep=4}, it suffices to notice that since $\sum_{\ell \geq 0} \Vert \E_0(X_{\ell} )\Vert_4< \infty$, we have the coboundary decomposition $S_n(f)=M_n(f) +r_0-r_0\circ \theta^n$ with $\| r_0 \|_4 < \infty$. So \eqref{conv1thgenep=4} follows directly from the fact that $(r_0-r_0\circ \theta^n)/n^{1/4}\to 0$ \as $\ $To prove \eqref{conv3thgenep=4}, we shall use Remark \ref{remarkdirectsense}. Therefore we need to show that
\begin{equation*}\label{p1conv3thgenep=4}
\sum_{k=1}^n (\E(d_k^2|\F_{k-1})-\E(d_k^2)) = O(n^{1/2} (\log \log n)^{1/2}) \qquad \as 
\end{equation*}
This condition follows directly from Theorem 12 of \cite{MPP} together with the fact that $(d_k)_{k \in {\mathbb N}}$ is a martingale differences sequence provided that
\begin{equation}\label{p2conv3thgenep=4}
\sum_{n \geq 2} \frac{(\log n)^3}{n^{2}}\|\E_0(M_n^2)-\E(M_n^2)\|_2^2 <\infty \, .
\end{equation}
According to the proof of Theorem 2.3 and Corollary 2.1 in \cite{DDM}, this will hold true provided that there exists $\gamma\in ]0,1]$ such that
\beq \label{cond1coralphaphi}
 \sum_{n>0} (\log n)^3 n^{ \frac{1}{\gamma} +\frac{1}{2}} \| \E_0(X_n) \|^{2}_4  < \infty \, ,
\eeq
and
\beq \label{cond2coralphaphi}
 \sum_{n>0}  (\log n)^3 n^{2\gamma }   \sup_{i \geq j \geq n }\|\E_0(X_iX_j ) - \E(X_iX_j)\|_{2}^{2} < \infty \, .
\eeq
Notice that $\| \E_0(X_n) \|_4 = \lim_{L \rightarrow \infty} \| \E_0(f_L (Y_n) - \E(f_L (Y_n) )) \|_4\le \liminf_{L\to \infty}
\sum_{k=1}^L|a_{k,L}|\, \|\E_0(f_{k,L}(Y_{n})-\E(f_{k,L}(Y_n))\|_4$. Next, by Lemma \ref{covaphi},  $\|\E_0(f_{k,L}(Y_{n})-\E(f_{k,L}(Y_n))\|_4\le
2 M (2\phi_{1, {\bf Y}}(n))^{3/4 }$. On the other hand, $$\|\E_0(X_iX_j ) - \E(X_iX_j)\|_{2} \leq \liminf_{L\to \infty}
\sum_{k=1}^L\sum_{\ell=1}^L |a_{k,L}||a_{\ell,L}| \, \|\E_0(f_{k,L}(Y_{i})f_{\ell,L}(Y_{j}) -\E(f_{k,L}(Y_{i})f_{\ell,L}(Y_{j})\|_2 \, ,$$ and by Lemma \ref{covaphi}
$$
\sup_{i \geq j \geq n } \|\E_0(f_{k,L}(Y_{i})f_{\ell,L}(Y_{j}) -\E(f_{k,L}(Y_{i})f_{\ell,L}(Y_{j})\|_2 \leq 16  M^2 (\phi_{2, {\bf Y}}(n))^{1/ 2} \, .
$$
Hence condition \eqref{condDDM} (for $p=4$) implies that \eqref{cond1coralphaphi} and \eqref{cond2coralphaphi} hold for $\gamma =1/{\sqrt 3}$. The proof of Theorem \ref{asipgenephi} is therefore complete. $\square$

\section{Proof of the results of Section \ref{sectionASIPUEM}} \label{sectionproofsectionASIPUEM}
\setcounter{equation}{0}

\subsection{Proof of Theorem \ref{ASmap1} on uniformly expanding  maps}

Item 1 follows directly from Theorem \ref{asipgenephi}. Indeed since $T$ is uniformly
expanding, it follows from Section 6.3 in \cite{DePr}
that
 \beq \label{majphi}
 \phi_{2, {\bf Y}}(n)=O(\rho^n) \ \text{ for some } \ \rho
\in (0, 1) \, .
\eeq

To prove Item 2 we proceed as follows. We start by the case $p\in ]2, 4[$. Since $f\in \Monm_p(M,\nu)$, we consider the function $f_j$ and the random variables $\bar X_{j,k}$ and $\widetilde X_{j,k}$ defined in the beginning of the proof of Theorem \ref{asipgenephi}. In addition we set 
\beq \label{deftribus}
{\mathcal F}_k = \sigma( Y_i, i \leq k) \, \text{ and } {\mathcal G}_k = \sigma( T^i, i \geq k) \, .
\eeq
As in the proof of Theorem \ref{asipgenephi}, we define a sequence of martingale differences, $(\bar d_{j,\ell})_{\ell \geq 1}$,  with respect to the non-decreasing sequence of  $\sigma$-algebras $({\mathcal F}_\ell)_{\ell\geq 1}$, as follows: $$
\bar d_{j,\ell} = \sum_{k \geq \ell}\big ( \E \left(\bar X_{j,k} | {\mathcal F}_{\ell}\right) -\E \left(\bar X_{j,k} | {\mathcal F}_{\ell-1}\right)  \big ) :=\sum_{k \geq \ell}P_{\ell} \left(\bar X_{j,k} \right) \, ,$$
and we recall that by Claim \ref{convlqdjk} and by \eqref{majphi},  the series $\sum_{k \geq 0}P_0\left(\bar X_{j,k}\right)$ converges in ${\mathbb L}^{\infty}$. Notice now that by the Markovian property of $(Y_i)_{i \geq 0}$, we have that
\beq \label{markovdjk}
\bar d_{j,\ell} = \sum_{k \geq \ell} \big ( \E \left(\bar X_{j,k} | Y_{\ell}\right) - \E \left(\bar X_{j,k} | Y_{\ell-1}\right) \big ) : = m(Y_{\ell}, Y_{\ell -1}) \, ,
\eeq
where $m( \cdot, \cdot)$ is a measurable function from ${\mathbb R}^2$ to ${\mathbb R}$. Define now
\beq \label{markovdjkMap}
 d^*_{j,\ell} =  m( T^{\ell -1}, T^{\ell}) \, .
\eeq
Notice then that $d^*_{j,\ell}$ is ${\mathcal G}_{\ell -1}$-measurable. Moreover, since on the probability space $([0, 1], \nu)$, the random
vector $(T, T^2, \ldots , T^n)$ is
distributed as $(Y_n,Y_{n-1}, \ldots, Y_1)$, we have that
$$
\|  \E (  d^*_{j,\ell} |{\mathcal G}_{\ell }) \|_{1, \nu} =  \| \E (  d^*_{j,\ell} |T^{\ell }) \|_{1, \nu} = \| \E (  \bar d_{j,\ell} |Y_{\ell-1 }) | \|_{1}  =0 \, ,
$$
it follows that $\E (  d^*_{j,\ell} |{\mathcal G}_{\ell })=0 $ $\nu$-a.s.

We define now some non stationary sequences $(X^*_{\ell})_{\ell \geq 1}$ and $( d^*_{\ell})_{\ell \geq 1}$ as follows:
\beq \label{def*1}
d^*_1:= d^*_{1,1}  \, , \   X^*_1:= \bar f_1 \circ T  - \nu (\bar f_1 \circ T ) \, ,
\eeq
and, for every $j\ge 0$ and every $\ell\in\{2^{j}+1,...,2^{j+1}\}$,
\beq \label{def**}  d^*_\ell:= d^*_{j,\ell}  \, , \
         X^*_\ell:=\bar f_j \circ T^{\ell}  - \nu (\bar f_j \circ T^{\ell}) \, .
\eeq
For every positive integer $n$, we then define
$$  M^*_n(f):=\sum_{\ell=1}^n  d^*_\ell\ \ \mbox{and}\ \  S^*_n(f):=\sum_{\ell=1}^n  X^*_{\ell} \, .$$
Therefore, Item 2 of Theorem \ref{ASmap1}  will follow if we can prove that 
\beq \label{conv1thgeneMap}
\sup_{1 \leq k \leq n } \big |\sum_{i=1}^k (f \circ T^i  - \nu (f) )  - \bar S^*_k(f)\big | = o (n^{1/p} ( \log n)^{1 - 2/p} ) \ \nuas \, ,
\eeq
\beq \label{conv2thgeneMap}
\sup_{1 \leq k \leq n } \big |\bar S^*_k(f)- \bar M^*_k(f)\big | = o (n^{1/p} ( \log n)^{1 - 2/p} ) \ \nuas \, ,
\eeq
and if, enlarging our probability space if necessary, there exists  a
    sequence $(Z^*_i)_{i \geq 0}$ of iid\ Gaussian random
    variables with mean zero and variance $\sigma^2$
such that
\begin{equation}\label{conv3thgeneMap}
\sup_{1 \leq k \leq n }\Big|\sum_{i=1}^{k} ( d^*_i
- Z^*_i) \Big| =o (n^{1/p} ( \log n)^{1 - 2/p} ) \ \nuas 
\end{equation}
According to the proof of Theorem \ref{asipgenephi}, \eqref{conv1thgeneMap} will hold if
 \begin{equation} \label{conv1thgenep1Map}
\sum_{j >0} 2^{-2j/p} j^{-2 + 4/p} \Big \|  \max_{1 \leq k \leq 2^j} \Big |  \sum_{\ell=1}^{k} \big ( (f- \bar f_j) \circ T^\ell  - \nu (f - \bar f_j)  \big ) \Big | \Big \|^2_{2, \nu} < \infty \, .
\end{equation}
But, since on the probability space $([0, 1], \nu)$, the random
vector $(T, T^2, \ldots , T^n)$ is
distributed as $(Y_n,Y_{n-1}, \ldots, Y_1)$, according to the inequality (4.1) in \cite{DeGoMe10}, we have
\[
\nu \Big( \max_{1 \leq k \leq
  n} \Big |\sum_{\ell=1}^{k} \big ( (f- \bar f_j) \circ T^\ell  - \nu (f - \bar f_j) ) \big )\Big|> x\Big)
  \leq
 \p \Big( 2 \max_{1 \leq k \leq
  n} \Big |\sum_{\ell=1}^{k} \widetilde X_{j, \ell}  \Big|> x \Big) \, .
\]
Therefore \eqref{conv1thgenep1Map} follows from Claim \ref{lmaphi1}.

We turn now to the proof of \eqref{conv2thgeneMap}. According to the proof of Theorem \ref{asipgenephi}, \eqref{conv2thgeneMap} will hold if
 \begin{equation} \label{conv2thgenep1Map}
\sum_{j >0} 2^{-4j/p} j^{-4 (1 -2/p)} \Big \|  \max_{1 \leq k \leq 2^j} \big |  \sum_{\ell=1}^{k} \big (  \bar f_j \circ T^{\ell}  - \nu ( \bar f_j)  - d^*_{j,\ell} \big ) \big | \Big \|^4_{4, \nu} < \infty \, .
\end{equation}
But, as before, since on the probability space $([0, 1], \nu)$, the random
vector $(T, T^2, \ldots , T^n)$ is
distributed as $(Y_n,Y_{n-1}, \ldots, Y_1)$,
\[\nu \Big( \max_{1 \leq k \leq 2^j} \big |  \sum_{\ell=1}^{k} \big (  \bar f_j \circ T^{\ell}  - \nu ( \bar f_j)  - d^*_{j,\ell}) \big ) \big | >  x\Big)
  \leq
 \p \Big( 2 \max_{1 \leq k \leq
  n} \Big |\sum_{\ell=1}^{k} (\bar X_{j, \ell}  - \bar d_{j,\ell}) \Big|> x \Big) \, .
\]
Therefore \eqref{conv2thgenep1Map} follows from Claim \ref{lmaphi2}.

To prove \eqref{conv3thgeneMap}, we shall proceed as for the proof of \eqref{conv3thgene} with the difference that Theorem \ref{shaoresult} is used instead of Theorem 2.1 in \cite{shao93}. So we have to prove that
$$
\sum_{i=1}^n\left({\mathbb E}((d_i^*)^2 | T^i)-{\mathbb E}((d_i^*)^2)\right)
     =o ( a_n )\ \ a.s. \, ,
$$
where $a_n = n^{2/p} ( \log n)^{1-4/p} $. But following the proof of \eqref{condshao2}, this will follow if we can prove that
\begin{equation} \label{condshao2p1Map}
\sum_{j >0} 2^{-4j/p} j^{-4 (1 -2/p)} \Big \|  \max_{1 \leq k \leq 2^j} \big |  \sum_{i=1}^{k} \left({\mathbb E}((d_{j,i}^*)^2 | T^i)-{\mathbb E}((d_{j,i}^*)^2)\right) \big |   \Big \|_{2, \nu}^2 < \infty \, .
\end{equation}
But using again the fact that,  on the probability space $([0, 1], \nu)$, the random
vector $(T, T^2, \ldots , T^n)$ is
distributed as $(Y_n,Y_{n-1}, \ldots, Y_1)$, we have that
\[ \nu \Big( \max_{1 \leq k \leq 2^j} \big |  \sum_{i=1}^{k} \left({\mathbb E}((d_{j,i}^*)^2 | T^i)-{\mathbb E}((d_{j,i}^*)^2)\right) \big | >  x \Big)
  \leq
\p \Big( 2 \max_{1 \leq k \leq 2^j} \big |  \sum_{i=1}^{k} \left({\mathbb E}_{i-1}({\bar d}^{\, 2}_{j,i})-{\mathbb E}({\bar d}^{\, 2}_{j,i}) \right) \big | >  x \Big)\, .
\]
Therefore \eqref{condshao2p1Map} follows by the fact that \eqref{condshao2p2} holds true. This ends the proof of Item 2 of Theorem \ref{ASmap1} when $p \in ]2 , 4[$.

We turn now to the proof of Item 2 when $p=4$. Let $X_i = f(Y_i) - \nu (f)$. According to the beginning of the proof of Theorem \ref{asipgenephi}, \eqref{majphi} implies that $\sum_{k \geq 1} \Vert \E (X_k |{\mathcal F}_0) \Vert_{4} < \infty $. We define a sequence $( d_{\ell})_{\ell \geq 1}$ of stationary martingale differences  with respect to the non-decreasing sequence of  $\sigma$-algebras $({\mathcal F}_\ell)_{\ell\geq 1}$ and that are in ${\mathbb L}^4 $ as follows: $$
 d_{\ell} = \sum_{k \geq \ell}P_\ell\left( X_{k}\right)  \, .$$By the Markovian property of $(Y_i)_{i \geq 0}$,
\begin{equation*} \label{markovdlp=4}
d_{\ell} = \sum_{k \geq \ell} \big ( \E \left(X_k | Y_{\ell}\right) - \E \left(X_k | Y_{\ell-1}\right) \big ): = m(Y_{\ell}, Y_{\ell -1}) \, ,\end{equation*}
where $m( \cdot, \cdot)$ is a measurable function from ${\mathbb R}^2$ to ${\mathbb R}$. We define now
\begin{equation*} \label{markovdlMapp=4}
d^*_{\ell} =  m( T^{\ell -1}, T^{\ell}) \, .
\end{equation*}
As before, notice that $d^*_{\ell}$ is ${\mathcal G}_{\ell -1}$-measurable and satisfies $\E (  d^*_{j,\ell} |{\mathcal G}_{\ell })=0 $ $\nu$-a.s. According to Corollary \ref{corrateASIP}, the result will then follow if we can prove that
\begin{equation} \label{p1item2}
\sum_{k=1}^n (\E( (d^*_k)^2|\G_{k})-\E((d^*_k)^2)) = O(n^{1/2} ( \log \log n)^{1/2}) \qquad  \nuas 
\end{equation}
But  $ \E( (d^*_k)^2|\G_{k}) = \E( (d^*_k)^2|T^k) := h(T^k) $ $\nuas$ where $h( \cdot )$ is a measurable function such that $\nu (  h^2) < \infty$. Let $ \widetilde h = h - \nu (h)$. Assume that we can prove that \begin{equation} \label{condMPP}
\sum_{n\geq 1}  (\log n)^3 \frac{\| \widetilde h +K \widetilde h+\cdots +K^{n-1} \widetilde h \|_{2, \nu}^2}{n^2} <\infty \, .
\end{equation}
This condition implies in particular that $\sum_{n\geq 1}  \frac{\| \widetilde h +K \widetilde h+\cdots +K^{n-1} \widetilde h \|_{2, \nu}}{n^{3/2}} <\infty$. By Lemma 2 of \cite{MPP}
(and its proof), it follows that
$\varphi(Y_0,Y_1): =\lim_n \frac{1}n \sum_{j=1}^n\sum_{k=0}^{j-1}(K^{k}\widetilde h(Y_1)
-K^k \widetilde h(Y_0))$ exists in $\LL^2$ and
$M_n (\varphi) :=\sum_{k=1}^n \varphi(Y_{k-1},Y_k)$ is a martingale with stationary
increments such that
$$
\| \widetilde h (Y_1)+\cdots +  \widetilde h (Y_n)-M_n (\varphi)\|_{2} \ll \sqrt n \sum_{k\ge n}
\frac{\| \widetilde h +K \widetilde h  +\cdots +K^{k-1} \widetilde h \|_{2, \nu }}{k^{3/2}}\, .
$$
Since  the random
vectors $(Id, T, T^2, \ldots , T^{n-1})$ and
 $(Y_n,Y_{n-1}, \ldots, Y_1)$ have same distribution, we can write  $M_n^* (\varphi):= \sum_{k=1}^n \varphi(T^k,T^{k-1})=\sum_{k=0}^{n-1}\varphi(T^{n-k},T^{n-k-1})$. Then $(M^*_n)$ is a sum associated to a stationary sequence of reverse martingale differences and 
\begin{equation} \label{majMPP}
\|\widetilde h +\cdots + \widetilde h \circ T^{n-1}-M_n^*(\varphi)\|_{2, \nu } \ll \sqrt n \sum_{k\ge n}
\frac{\| \widetilde h +K \widetilde h +\cdots +K^{k-1} \widetilde h \|_{2, \nu}}{k^{3/2}}\, .
\end{equation}
Therefore, by Corollary 4.2 of \cite{Cuny} with $b(n)=\log n$, if
\begin{equation} \label{condMPPbis}
\sum_n \frac{\log n}{n^2} \|\widetilde h +\cdots + \widetilde h \circ T^{n-1} -M_n^*(\varphi)\|_{2,\nu }^2< \infty \, ,
\end{equation}
then
$$
\frac{\widetilde h +\cdots + \widetilde h \circ T^{n-1} -M_n^*(\varphi)}{\sqrt{n\log\log n}}\to 0 \qquad \mbox{$\nu$-a.s.}
$$
Using \eqref{majMPP}, it is easy to see that \eqref{condMPPbis} holds as soon as \eqref{condMPP} is satisfied. Using then Corollary \ref{corFLIL}  to observe  that $M_n^*(\varphi)  = O ( \sqrt{n\log\log n} ) $ $ \nuas$, we conclude that \eqref{p1item2} (and then Item 2 when $p=4$) holds as soon as \eqref{condMPP} does. Notice now that \eqref{condMPP} can be rewritten as
$$
\sum_{n\geq 1}  (\log n)^3 \frac{\| \sum_{k=1}^n (\E( d_k^2|\F_{0})-\E(d_k^2))  \|_{2}^2}{n^2} <\infty \, ,
$$
that is exactly condition \eqref{p2conv3thgenep=4}.  $\square$

\subsection{Proof of Theorem \ref{ASmapMC}} The proof follows directly by analyzing the proof of Theorem \ref{ASmap1} when $p=4$. Indeed, the proof reveals that if the conditions \eqref{cond1coralphaphi} and \eqref{cond2coralphaphi} hold for $X_i = f(Y_i) - \nu (f)$, then the strong approximation principle holds for both $\sum_{i=1}^n (f(Y_i) - \nu (f))$ and $\sum_{i=1}^n (f \circ T^i - \nu (f))$ with rate $O(n^{1/4}(\log n)^{1/2} (\log n \log n)^{1/4} )$. The condition \eqref{cond1ASmapMC} (resp. \eqref{cond2ASmapMC}) is exactly the condition \eqref{cond1coralphaphi} (resp. \eqref{cond2coralphaphi}) rewritten with the help of the transition operator $K$. $\square$

\subsection{Proof of Theorem \ref{ASmapMCconcave} on uniformly expanding  maps} It suffices to verify the assumptions of Theorem \ref{ASmapMC}. Using the first part of Condition \eqref{cond1PFO} and the fact that $ | f(x) - f(y) | \leq c( |x-y|) $ where $c$ is a concave and non-decreasing function,  it follows from Lemma 17 in \cite{DMPU} that $\| K^n (f) - \nu (f) \|_{\infty,\nu}  \leq c(C \rho^n) $. Therefore, \eqref{cond1ASmapMC} is satisfied with $\gamma=1/{\sqrt 3}$ as soon as \eqref{intc} is.  To verify now the condition \eqref{cond2ASmapMC} of Theorem \ref{ASmapMC}, we shall use similar arguments as those developed in the proof of Corollary 3.12 in \cite{DDM}. From Section 7 in \cite{DP07}, we know that for $i$ and $j$ positive integers, there exists $(Y_i^*, Y_j^*)$ distributed as $(Y_i, Y_j)$ and independent of $Y_0$ such that
$$
  \frac12 \big \|  {\mathbb E}(|Y_i-Y_i^*| | Y_0) + {\mathbb E}(|Y_j-Y_j^*| | Y_0) \big \|_{\infty}=
  \sup_{h \in \Lambda_1({\mathbb R}^2)}\ \big \|{\mathbb E}(h(Y_i, Y_j)|Y_0)-{\mathbb E}(h(Y_i, Y_j)) \big \|_{\infty} \, .
$$
Notice now that for $i \geq j \geq 0$, by using the second part of Condition \eqref{cond1PFO},
$$
 \sup_{h \in \Lambda_1({\mathbb R}^2)}\big \|{\mathbb E}(h(Y_i, Y_j)|Y_0)-{\mathbb E}(h(Y_i, Y_j)) \big \|_{\infty} =  \sup_{h \in \Lambda_1({\mathbb R}^2)} \Vert K^j \circ Q_{i-j}  (h)
- \nu \big ( Q_{i-j}  (h) \big ) \Vert_{\infty, \nu}  \leq C \rho^j\, .
$$
On an other hand, we clearly have that
\begin{multline*}
\| K^j (f K^{i-j}(f)) - \nu (f K^{i-j}(f) ) \|_{\infty, \nu } = \big \|{\mathbb E}(f(Y_i) f(Y_j)-f(Y^*_i) f( Y^*_j)|Y_0) \big \|_{\infty}
\\
= \big \|{\mathbb E}( ( f(Y_i) - f(Y_i^*)) f(Y_j) |Y_0) -{\mathbb E}( f(Y_i^*)  ( f(Y_j^*) - f(Y_j))| Y_0)  \big \|_{\infty}  \, .
\end{multline*}
Since $f$ is continuous on a compact set, there exists a positive constant $R$ such that $\| f \|_{\infty} \leq R$. Hence,
$$
\| K^j (f K^{i-j}(f)) - \nu (f K^{i-j}(f) ) \|_{\infty, \nu } \leq R  \big \|{\mathbb E}( ( c ( | Y_i - Y_i^*| ) |Y_0) \big \|_{\infty}  + R \big \| {\mathbb E}( ( c ( | Y_j - Y_j^*| ) |Y_0)   \big \|_{\infty}  \, .
$$
Since $c$ is concave and and non-decreasing function, it follows that
\begin{align*}
\| K^j (f K^{i-j}(f))  - \nu (f K^{i-j}(f) ) \|_{\infty, \nu } & \leq R \big \|c \big ( {\mathbb E}( | Y_i - Y_i^*|  |Y_0) \big ) \big \| _{\infty} + R \big \|c \big ( {\mathbb E}( | Y_j - Y_j^*|  |Y_0) \big ) \big \| _{\infty}  \\
& \leq R \, c \big (  \big \|{\mathbb E}( | Y_i - Y_i^*|  |Y_0) \ \big \| _{\infty} \big )  + R c  \big ( \big \| {\mathbb E}( | Y_j - Y_j^*|  |Y_0)  \big \| _{\infty}   \big )\\
& \leq 2 R \, c \Big (  \frac 12 \big \|{\mathbb E}( | Y_i - Y_i^*|  |Y_0) \ \big \| _{\infty}  + \frac 12 \big \| {\mathbb E}( | Y_j - Y_j^*|  |Y_0)  \big \| _{\infty}   \Big ) \, .
\end{align*}
So overall,
$$
\sup_{i \geq j \geq n }  \| K^j (f K^{i-j}(f)) - \nu (f K^{i-j}(f) ) \|_{\infty, \nu } \leq   2 R \,  c ( C \rho^n) \, ,
$$
implying that \eqref{cond2ASmapMC} is satisfied with $\gamma=1/{\sqrt 3}$ as soon as \eqref{intc} is. This ends the proof of the theorem. $\square$

\section{Proofs of the reverse martingale's results} \label{sectproofmartingale}

 \setcounter{equation}{0}
 We start by recalling the following estimate of Hanson and Russo \cite[Theorem 3.2A]{HR}

\begin{lma}
Let $(B_t)_{t\ge 0}$ be a standard Brownian motion. Then
\beq\label{HR}
\lim_{a\rightarrow \infty} \sup_{t\geq 0} \sup_{0\leq s\leq a} \frac{|B_{t+s}-B_t|}{(2a [\log ((t+a)/a))+\log \log a])^{1/2}}
=1 \qquad \as
\eeq
\end{lma}
\medskip

We also recall the following convergence result for reverse martingales. \begin{lma}\label{reverse}
Let $(\xi_n)_{n \geq 1}$ be a sequence of variables in ${\mathbb L}^p$, $1\le p\le 2$, adapted
to a non-increasing filtration $(\G_n)_{n \geq 1}$.
 Assume that $\E(\xi_n|\G_{n+1})=0$
and $\sum_{ n\ge 1} \E(|\xi_n|^p)<\infty$.  Then $\sum_{n \geq 1} \xi_n$ converges
\as and in ${\mathbb L}^p$.
\end{lma}

\noindent{\bf Proof of Lemma \ref{reverse}.} The result is clear when $p=1$, hence we assume $p>1$.
Notice that for every $n>1$, $(\sum_{k=n-l}^{n}\xi_k)_{0\le l \le n-1}$ is a
$(\G_{n-l})_{0\le l \le n-1}$-martingale. Hence, by Burkholder inequality and
using that $x\mapsto |x|^{p/2}$ is subadditive, it follows that the exists a positive constant $C_p$ such that,
for every $1\le r < n$,
\begin{gather}\label{inmax}
\E \Big ( \max_{r\le m \le n}|\sum_{k=m}^n \xi_k|^p \Big )
\le C_p \E \Big ( \Big(\sum_{k=r}^n\xi_k^2\Big )^{p/2}\Big ) \le C_p
\sum_{k=r}^n \E(|\xi_k|^p) \, .
\end{gather}
So, $(Z_n):=(\sum_{k=1}^n \xi_k)$ is Cauchy in $\LL^p$, hence converges
in $\LL^p$, say to $Z$. Moreover, letting $n\to \infty$
in \eqref{inmax}, we see that  for every $r\ge 1$
\begin{equation*}
\E(\max_{m\ge r} |Z_m-Z|^p)\le C_p \sum_{k\ge r} \E(|\xi_k|^p)\, ,
\end{equation*}
which implies the desired result. \hfill $\square$

\subsection{Proof of Proposition \ref{huggins}.}
The $\LL^2$ and a.s. convergence of $\sum_{k\ge 1}\xi_k$ follows from Lemma \ref{reverse}. By Theorem 2 of Scott and Huggins \cite{HS}, enlarging our probability space if necessary, there exists a Brownian motion $(B_t)_{t\ge 0}$,
a non-increasing filtration $(\h_n)_{n \in {\mathbb N}}$ and a non-increasing process
$(\tau_n)_{n \in {\mathbb N}}$ adapted to $(\h_n)_{n \in {\mathbb N}}$, such that 
$$
R_n=B_{\tau_n} \qquad \as \,
$$
Moreover, writing $t_n:=\tau_n-\tau_{n+1}\ge 0$ $\P$-a.s., we have
\begin{gather}
\label{id}\E(t_n|\h_{n+1})=\E(\xi_n^2|\G_{n+1}) \qquad \as \, ,\\
\label{moment}\E(t_n^{p/2}|\h_{n+1})\le C_p \E(|\xi_n|^p|\G_{n+1}) \qquad \as
\quad \mbox{for every $p>1$ \, .}
\end{gather}

\medskip

Hence, using \eqref{id} twice,
$$\tau_n-\E(\tau_n)= \tau_n-\delta_n^2=\sum_{k\ge n} \big (t_k-\E(t_k|\h_{k+1}) \big )  +V_n^2-\delta_n^2
\qquad \as
$$
But it follows from \eqref{cond2} and \eqref{moment} that $\sum_{n\ge 1} \alpha_n^{-\nu} \E(t_n^{\nu}
)<\infty $ which implies, by Lemma \ref{reverse}, that $\sum_{k\ge 1}\alpha_k^{-1}(t_k-\E(t_k|\h_{k+1}))$ converges $\as$.
 Then, by an analogue to the Kronecker lemma
(see e.g. Heyde \cite[Lemma 1]{Heyde}),
$\sum_{k\ge n} (t_k-\E(t_k|\h_{k+1}))
=o(\alpha_n)$ $\as$ Together with \eqref{cond1}, this implies in particular that $\tau_n-\delta_n^2=o(\alpha_n)$ $\as$

\smallskip

For every $t>0$ define $\widetilde B_t=tB_{1/t}$, and $\widetilde B_0=0$. It is well-known that $(\widetilde B_t)_{t \ge 0}$ is a
standard Brownian motion. We have $B_{\tau_n}-B_{\delta_n^2}= \tau_n(\widetilde B_{1/\tau_n}-\widetilde B_{1/\delta_n^2})
+ (\tau_n-\delta_n^2)\widetilde B_{1/\delta_n^2}$. By the law of the iterated logarithm
for $(\widetilde B_t)_{t\ge 0}$ (or using that the supremum in \eqref{HR} is greater than
what we have for $t=0$ and $s=1/\delta_n^2$), we see that $\widetilde B_{1/\delta_n^2} =O(\delta_n^{-1}(\log \log
(1/\delta_n))^{1/2})$ \as Hence since $\alpha_n =O(\delta_n^{2})$, $(\tau_n-\delta_n^2)\widetilde B_{1/\delta_n^2}=o((\alpha_n \log \log
(1/\alpha_n))^{1/2})$ \as$\, $\\ Let us deal now with $\tau_n(\widetilde B_{1/\tau_n}-\widetilde B_{1/\delta_n^2})$. With this aim, we shall use \eqref{HR}. Since $\alpha_n =O(\delta_n^{2})$ and $\tau_n-\delta_n^2=o(\alpha_n)$ $\as$, we have $|1/\tau_n-1/\delta_n^2|=o(1/\alpha_n)$ \as
\quad  Define $u_n:=\alpha_n/(\tau_n\delta_n^2)$, $\varepsilon_n:=
\max(|\delta_n^2-\tau_n|/\alpha_n,u_n^{-1/2})$, $s_n:=|1/\tau_n-1/\delta_n^2|$,  $a_n:=\varepsilon_n u_n$  and $v_n:=\min
(1/\delta_n^2,1/\tau_n)$. Notice that $a_n\to \infty$,  $\varepsilon_n\to 0$,
$v_n+s_n=\max (1/\delta_n^2,1/\tau_n)$ and  $|\widetilde B_{1/\tau_n}-\widetilde B_{1/\delta_n^2}|=
|\widetilde B_{v_n+s_n}-\widetilde B_{v_n}|$. By \eqref{HR}, we have
\begin{gather*}
\frac{|\widetilde B_{v_n+s_n}-\widetilde B_{v_n}|}{(2a_n [\log
 ((v_n+a_n)/a_n)+\log \log a_n])^{1/2}}\\
\le \sup_{t\ge 0} \sup_{0\le s\le a_n} \frac{|\widetilde B_{t+s}-
\widetilde B_t|}{(2a_n[\log ((t+a_n)/a_n))+\log \log a_n])^{1/2}}
\to 1 \qquad \as
\end{gather*}
In particular, we have $|\widetilde B_{1/\tau_n}-\widetilde B_{1/\delta_n^2}|
=O \big( [\varepsilon_n u_n (|\log (\delta_n^2/(\alpha_n \varepsilon_n))|+
\log \log (u_n\varepsilon_n))]^{1/2} \big)$ \as $\, $
Then, using that $ |\log (\delta_n^2/(\alpha_n \varepsilon_n))|
\le |\log (\delta_n^2/\alpha_n)|+|\log \varepsilon_n|$ and that
$\varepsilon_nu_n\log\log (\varepsilon_nu_n)=o(u_n\log\log u_n)$, we obtain
\beq  \label{tildeBtaun}
\tau_n(\widetilde B_{1/\tau_n}-\widetilde B_{1/\delta_n^2})
=o\big( [\alpha_n  (|\log (\delta_n^2/\alpha_n)| +
\log \log (\alpha_n/\delta_n^4)]^{1/2}\big) \qquad \as\, ,
\eeq
which proves the result, since $1/\delta_n^4=O(1/\alpha_n^2)$. \hfill $\square$
\begin{rmk} \label{remlil}
It follows from the proof that the assumption \eqref{cond2} may be replaced by $\sum_{k\ge n} (t_k-\E(t_k|\h_{k+1}))
=o(\alpha_n)$ \as
\end{rmk}

\subsection{Proof of Theorem \ref{shaoresult}.}
Define $\xi_n:=X_n/\sigma_n^2$. Then, since $\E(\xi_k^2) = (\sigma^2_k -\sigma^2_{k-1})\sigma_k^{-4}$, by  comparing sums and integrals, it follows that  $\sum_{k \geq 1} \E(\xi_k^2)<\infty$.
Using the notations $V_n^2 = \sum_{k\geq n} (\E(\xi_k^2|\G_{k+1})$ and $\delta_n^2 =  \sum_{k\geq n} \E(\xi_k^2)$,  and writing
$T_n:=\sum_{k=1}^n (\E(X_k^2|\G_{k+1})-\E(X_k^2))$, we have
$$
V_n^2-\delta_n^2= \sum_{k\ge n}\frac{T_k-T_{k-1}}{\sigma_k^4}= \sum_{k\ge n}T_k \Big (
\frac{1}{\sigma_k^4}-\frac{1}{\sigma_{k+1}^4} \Big ) - \frac{T_{n-1}}{\sigma_n^4 }\, .
$$
Using \eqref{cond1'} and that $(\sigma_n)_{n \in {\mathbb N}}$ and $(a_n/\sigma_n^2)_{n \in {\mathbb N}}$ are respectively
non-decreasing and non-increasing, we obtain
$$
|V_n^2-\delta_n^2|=o\Big (\frac{a_n}{\sigma_n^2}\Big )\sum_{k\ge n}\Big (\frac{1}{\sigma_k^2}-\frac{1}{\sigma_{k+1}^2}\Big )+o\Big (\frac{a_n}{\sigma_n^4}\Big  )=o\Big (\frac{a_n}{\sigma_n^4}\Big )\, .
$$
We want to apply Proposition \ref{huggins} to $(\xi_n)$ with $\alpha_n:=a_n/\sigma_n^4$. Using \eqref{cond2'}, we have
$$
\sum_{i\ge 1} \alpha_i^{-\nu} \E(|\xi_i|^{2\nu})
=\sum_{n\ge 1} a_i^{-\nu} \E(|X_i|^{2\nu})<\infty \, ,
$$
hence condition \eqref{cond2} holds. It remains to prove that
$\alpha_n=O(\delta_n^2)$ and that $\alpha_n /\delta_n^4\to \infty$. With this aim, we first notice that
\begin{gather*}
\frac{1}{\sigma_{n-1}^2}-\delta_n^2 =\sum_{k\ge n} \int_{\sigma_{k-1}^2}^{\sigma_k^2}
\Big (\frac1{x^2}-\frac1{\sigma_k^4}\Big )dx \, .
\end{gather*}
Hence, using that $\sup_n \E(X_n^2)<\infty$, it follows that
$\sigma_n=O(\sigma_{n-1})$ and
\begin{gather*}
0\le \frac{1}{\sigma_{n-1}^2}-\delta_n^2\le \sum_{k\ge n}
\E(X_k^2) \Big (\frac1{\sigma_{k-1}^4}-\frac1{\sigma_k^4} \Big )
=O \Big ( \sum_{k\ge n}
\frac{\sigma_k^2-\sigma_{k-1}^2}{\sigma_k^6} \Big ) =
O \Big( \frac{\delta_n^2}{\sigma_{n}^2} \Big )\, .
\end{gather*}
In particular, since $\delta_n^2=O(\sigma_{n-1}^{-2}) =O(\sigma_{n}^{-2}) $ and $|\sigma_n^{-2} -\sigma_{n-1}^{-2} | = O( \sigma_n^{-4})$, we have
\begin{equation}\label{sigma_n}
\Big |\frac{1}{\sigma_{n}^2}-\delta_n^2  \Big|=O\Big (\frac{1}
{\sigma_n^4} \Big )\, .
\end{equation}
Since $a_n \sigma_n^{-2}$ is non-increasing, \eqref{sigma_n} implies that $\alpha_n=O(\delta_n^2)$. In addition since $a_n$ is tending to infinity,  \eqref{sigma_n} entails also that  $\alpha_n /\delta_n^4\to \infty$.

\medskip

By Proposition \ref{huggins}, enlarging our probability
space if necessary, there exists a standard Brownian motion
$(B_t)_{t\ge 0}$, such that \eqref{raterenv} holds with $\delta_n^2=\sum_{k\ge n}
(\sigma_k^2-\sigma_{k-1}^2)\sigma_k^{-4}$.  Now, for every $t>0$ define $\widetilde B_t=tB_{1/t}$, and $\widetilde B_0=0$ (recall that $(\widetilde B_t)_{t \ge 0}$ is a
standard Brownian motion). Notice that $$B_{1/\sigma_n^2} -B_{\delta_n^2}= \sigma_n^{-2}(\widetilde B_{\sigma_n^2}-\widetilde B_{1/\delta_n^2})
+ (\sigma_n^{-2}-\delta_n^2) \widetilde B_{1/\delta_n^2}\, .$$ By \eqref{sigma_n} and the law of the iterated logarithm
for $(\widetilde B_t)_{t\ge 0}$, we derive that $$
(\sigma_n^{-2}-\delta_n^2)\widetilde B_{1/\delta_n^2 } = O ( \sigma_n^{-2}(\log \log
(\sigma_n))^{1/2}) = o   \bigg(\frac{ \big (a_n (\log \log a_n )  \big )^{1/2}}{\sigma_n^2}\bigg) \ \as  $$
To deal now with $\sigma_n^{-2}(\widetilde B_{\sigma_n^2}-\widetilde B_{1/\delta_n^2})$, we use the same arguments as the ones used to derive \eqref{tildeBtaun} (with $\sigma_n^{-2} $ replacing $\tau_n$). Hence we infer that
$$
\sigma_n^{-2}(\widetilde B_{\sigma_n^2}-\widetilde B_{1/\delta_n^2}) = o \bigg(\frac{ \big (a_n(|\log (\sigma_n^2/a_n)|+\log \log a_n  ) \big )^{1/2}}{\sigma_n^2}\bigg) \qquad \mbox{$\P$-a.s.}
$$
So, overall, it follows that
\beq \label{resinterRn}
|R_n-B_{1/\sigma_n^2}|= o \bigg(\frac{ \big (a_n(|\log (\sigma_n^2/a_n)|+\log \log a_n  ) \big )^{1/2}}{\sigma_n^2}\bigg) \qquad \mbox{$\P$-a.s.} \, ,
\eeq
where $R_n= \sum_{k \geq n}X_n/\sigma_n^2$.

Write $\widetilde Z_n:=\sigma_n^2(B_{1/\sigma_n^2}-B_{1/\sigma_{n+1}^2})$. By independence of the increments,
$(\widetilde Z_n)$ is a sequence of independent centered Gaussian variables.
Notice that, by stationarity of the increments
$\E(\widetilde Z_n^2)=\sigma_n^2\E(X_n^2)/\sigma_{n+1}^2$.

\smallskip

We have
\begin{align*}
\sum_{k=1}^n X_k-\sum_{k=1}^n \widetilde Z_k&=\sum_{k=1}^n \sigma_k^2(
(R_k-B_{1/\sigma_k^2})-(R_{k+1}-B_{1/\sigma_{k+1}^2}) )\\
&=\sum_{k=2}^n (R_k-B_{1/\sigma_k^2})(\sigma_k^2-\sigma_{k-1}^2)
+\sigma_1^2(R_1-B_{1/\sigma_1^2})
-\sigma_n^2(R_{n+1}-B_{1/\sigma_{n+1}^2})\, .
\end{align*}
Using that $(\sigma_n)$, $(\sigma_n^2/a_n)$, $(a_n/\sigma_n)$ and $(a_n)$ are
non-decreasing, and taking into account \eqref{resinterRn}, we deduce that
\begin{align*}
\sum_{k=1}^n X_k-\sum_{k=1}^n \widetilde Z_k
& =o \bigg(\Big(\frac{a_n(|\log (\sigma_n^2/a_n)| +\log \log a_n  ) }{\sigma_n }\Big)^{1/2}\bigg) \sum_{k=2}^n \frac{\sigma_k^2-\sigma_{k-1}^2}{(\sigma_k^2
)^{3/4}}\\
& =o \bigg(\big(a_n(|\log (\sigma_n^2/a_n)|+\log \log a_n
\big)^{1/2}\bigg)\, ,
\end{align*}
where we used that $\sum_{k=2}^n (\sigma_k^2-\sigma_{k-1}^2) (\sigma_k^2
)^{-3/4}=O(\int_0^{\sigma_n^2}dx/x^{3/4})$.

\smallskip

Finally, define $Z_n:=\widetilde Z_n\sigma_{n+1}/\sigma_n $. Notice that
$|Z_n-\widetilde Z_n|\le C|\widetilde Z_n|/\sigma_n$ for some $C>0$. Hence
$\sum_{n \geq 1} \E((Z_n-\widetilde Z_n)^2)/a_n\le \sum_{n\ge 1}(\sigma_{n}^2-\sigma_{n-1}^2)/\sigma_n^3<\infty$. So, by the Kolmogorov theorem
(see also Lemma \ref{reverse}),
$\sum_{n \geq 1} (Z_n-\widetilde Z_n)/\sqrt{a_n}$
converges  $\P$-a.s., and \eqref{rate} follows from
the Kronecker lemma. \hfill $\square$

\subsection{Proof of Corollary \ref{corFLIL}.} Assume that $\E (X_1^2) \neq 0$, otherwise there is nothing to prove.  We start by the proof of Corollary \ref{corFLIL}. Notice first that by stationarity and Fubini theorem,
\begin{gather*}
\sum_{n\ge 1} \frac{\E ( |X_n|{\bf 1}_{\{|X_n| >\sqrt n\}})}{\sqrt n}
= \E \Big  (|X_1| \sum_{1\le n <X_1^2}\frac{1}{\sqrt n} \Big  )
\le C \E(X_1^2)<\infty \, .
\end{gather*}
Hence,
\begin{equation}\label{neweq1}
\sum_{n\ge 1} n^{-1/2} |X_n|{\bf 1}_{\{|X_n| >\sqrt n\}}
<\infty  \quad \as  \,
\mbox{ and } \,
\sum_{n\ge 1} n^{-1/2} \E(|X_n|{\bf 1}_{\{|X_n| >\sqrt n\}}|\G_{n+1})
<\infty  \quad \as
\end{equation}
and by the Kronecker lemma,
\begin{equation}\label{neweq2}
\sum_{k= 1}^n |X_k|{\bf 1}_{\{|X_k| >\sqrt k\}}=o(\sqrt n)\quad \as
 \, \mbox{ and }  \,
\sum_{k= 1}^n \E(|X_k|{\bf 1}_{\{|X_k| >\sqrt k\}}|\G_{k+1})=o(\sqrt n)\quad \as
\end{equation}

\medskip

Define $Y_n:=X_n{\bf 1}_{\{|X_n|\le \sqrt n\}}-
\E(X_n{\bf 1}_{\{|X_n|\le \sqrt n\}}|\G_{n+1})$. Then, by the above,
using that $\E(X_n|\G_{n+1})=0$ a.s., we see that it suffices to
prove \eqref{FLIL} with $(X_n)$ replaced with $(Y_n)$.

\medskip

We want to apply Theorem \ref{shaoresult} to $(Y_n)$ with $a_n=\sigma^2_n=n$.
We have to prove conditions \eqref{cond1'} and \eqref{cond2'}. Let us prove \eqref{cond1'}. Clearly,
$(\E (Y_1^2)+\cdots + \E(Y_n^2))/n\to  \E(X_1^2)$. Hence, we
only need to prove that
\begin{equation}\label{EYn}
(\E (Y_1^2|\G_2)+\cdots + \E(Y_n^2|\G_{n+1}))/n
\to  \E(X_1^2)\qquad  \as
\end{equation}
We first prove that 
\begin{equation}\label{EYnbis}
(\E (Y_1^2|\G_2)+\cdots + \E(Y_n^2|\G_{n+1}))-(Y_1^2+\cdots + Y_n^2)
=o(n)\qquad  \as
\end{equation}By Kronecker lemma, this will follow from the convergence of the series
$\sum_n (\E(Y_n^2|\G_{n+1})-Y_n^2)/n$. By Lemma \ref{reverse}, this last convergence will hold true provided that $\sum_n \E(Y_n^4)/n^2 <\infty$.  
But, by stationarity and Fubini theorem, we have
\begin{gather}\label{fub}
\sum_{n\ge 1}\frac{\E(Y_n^4) }{n^2}
\le  16 \, \E \Big (X_1^4 \sum_{n\ge X_1^2} \frac{1}{n^2} \Big )\le C\E(X_1^2)<\infty \, .
\end{gather}
Therefore \eqref{EYnbis} is proved. 
Now, by the ergodic theorem we have
$$
\limsup_n \frac{\sum_{k=1}^n X_k^2{\bf 1}_{\{|X_k|\le \sqrt k\}}}n \le
\lim_n \frac{\sum_{k=1}^n X_k^2}n=\E(X_1^2)
\qquad \as\, ,
$$
and for any $A$ fixed,
$$\liminf_n\frac{\sum_{k=1}^n X_k^2{\bf 1}_{\{|X_k|\le \sqrt k\}}}n \ge
\lim_n \frac{\sum_{k=1}^n X_k^2{\bf 1}_{\{|X_k|\le A\}}}n=\E(X_1^2{\bf 1}_{\{|X_1|\le A\}})
\qquad \as
$$
Letting $A\to \infty$, we see that the $\liminf$ and the $\limsup$ above
are equal to $\E(X_1^2)$. Hence
\begin{equation}\label{neweq3}
\frac{\sum_{k=1}^n X_k^2{\bf 1}_{\{|X_k|\le \sqrt k\}}}n \to \E(X_1^2) \qquad \as
\end{equation}
On an other hand using the fact that  $\E(X_k|\G_{k+1}) = 0$ a.s. together with \eqref{neweq2}, we get that
\begin{equation}\label{neweq4}
n^{-1}\sum_{k=1}^n  (\E(X_k{\bf 1}_{\{|X_k| \leq \sqrt k\}}|\G_{k+1}))^2 \le n^{-1/2}\sum_{k=1}^n  \E( |X_k |{\bf 1}_{\{|X_k| > \sqrt k\}}|\G_{k+1})
=o(1) \qquad \as
\end{equation}
Combining \eqref{neweq3}, \eqref{neweq4} and \eqref{EYnbis}, we see that \eqref{EYn} holds, which proves \eqref{cond1'}.

\medskip

The fact that \eqref{cond2'} holds with $\nu=2$ follows from
\eqref{fub}.
By Theorem \ref{shaoresult}, there exists
a sequence of independent centered Gaussian variables $(\widetilde Z_n)_{n \geq 1}$ such that
$\E(\widetilde Z_n^2)=\E(Y_n^2)=\E(X_1^2)+o(1)$ and $Y_1+ \cdots Y_n-(\widetilde Z_1+\cdots +\widetilde Z_n)=o(\sqrt{n\log \log n})$ \as $\ $Let $(\delta_k)_{k \geq
1}$ be a sequence of iid\ Gaussian random variables with mean zero 
and variance $\E (X_1^2)$, independent of the sequence $(\widetilde Z_{n})_{n\geq
1}$. We now construct a sequence $(Z_{n})_{n \geq 1}$ as follows. If $\E(\widetilde Z_n^2)=0$, then $Z_n=\delta_n$, else $Z_n=
c_n \widetilde Z_{n}$ where $c_n = \sqrt{\frac{\E(X_1^2)}{\E(\widetilde Z_n^2)}}$.  By construction, the $Z_n$'s are iid\ Gaussian random variables with
mean zero and variance $\E(X_1^2)$. Write $G_n:=Z_n-\widetilde Z_n$ and  $v^2_n:=\sum_{k=1}^n \E (G_{k}^2)$. By L\'evy's inequality (see for instance Proposition 2.3 in \cite{LT}),
\beq \label{inelevy}
  {\mathbb P}\Big
  ( \max_{1 \leq k \leq 2^r} \Big| \sum_{i=1}^k G_i \Big|>  x\Big) \leq  2 \exp \Big(-\frac{ x^2}{2  v^2_{2^r}}\Big)\, .
\eeq
Hence taking $x= 2 v_{2^r}
  (\log \log 2^r)^{1/2}$, we get that
  \begin{equation*}
   \sum_{r\ge 0} {\mathbb P}\Big
  ( \max_{1 \leq k \leq 2^r} \Big| \sum_{i=1}^k G_i \Big|>  2 v_{2^r}
  (\log \log 2^r)^{1/2} \Big) < \infty \, .
  \end{equation*}
Therefore $\sup_{1 \leq k \leq 2^r} \big| \sum_{i=1}^k G_i \big|= O \big ( v_{2^r} (\log \log 2^r)^{1/2} \big)$ almost surely.  
\as $ \ $ This ends the proof of Corollary \ref{corFLIL} since $v^2_n=o(n)$. \hfill $\square$

\medskip

\subsection{Proof of Corollary \ref{corrateASIP}.}
Define $Y_n:=X_n{\bf 1}_{\{|X_n|\le  n^{1/p}\}}-
\E(X_n{\bf 1}_{\{|X_n|\le n^{1/p}\}}|\G_{n+1})$. Since $\E (X_n | \G_{n+1}) = 0$ a.s.,
$$
\sum_{k\geq 1} \frac{\E |X_k-Y_k|}{k^{1/p}} \leq 2\sum_{k\geq 1} \frac{\E (  |X_k|{\bf 1}_{\{|X_n|\le  k^{1/p}\}})}{k^{1/p}} \, .
$$
Hence by stationary and Fubini theorem, $\sum_{k\geq 1} k^{-1/p}\E |X_k-Y_k| < \infty $, implying via the Kronecker lemma that
$$
\sum_{k=1}^n |X_k-Y_k| =o(n^{1/p}) \qquad \as
$$
Let us prove now that $(Y_n)_{n \geq 1}$ satisfies the conditions of Theorem \ref{shaoresult} with $a_n = n^{2/p} b(n)$ and $\sigma_n^2=n$. With this aim, we first notice that since $\E (X_n | \G_{n+1}) = 0$ a.s.,
$$
\E(X_k^2|\G_{k+1})-\E(X_k^2) -\E(Y_k^2|\G_{k+1})+\E(Y_k^2) = \E(X_k^2{\bf 1}_{\{|X_k|> k^{1/p}\}}|\G_{k+1})-
\E(X_k^2{\bf 1}_{\{|X_k|> k^{1/p}\}}) \, .$$
Since by stationarity and Fubini theorem,  $\sum_{k\geq 1} k^{-2/p}\E(X_k^2{\bf 1}_{\{|X_k|> k^{1/p}\}}) < \infty $, we conclude via the Kronecker lemma that
\begin{gather*}
\sum_{k=1}^n |\E(X_k^2|\G_{k+1})-\E(X_k^2) -\E(Y_k^2|\G_{k+1})+\E(Y_k^2)|
=o(n^{2/p}) \qquad \as
\end{gather*}
Together with condition \eqref{condcarrecondi}, this implies that
\begin{equation*}
\sum_{k=1}^n (\E(Y_k^2|\G_{k+1})-\E(Y_k^2))
=o(n^{2/p}b(n)) \qquad \as
\end{equation*}
Notice now that by stationarity and Fubini theorem,
\begin{gather*}
\sum_{n\ge 1}\frac{\E(Y_n^4) }{n^{4/p}}
\le  16 \, \E \Big (X_1^4 \sum_{n\ge |X_1|^p} \frac{1}{n^{4/p}} \Big )\le C_p\E(|X_1|^p)<\infty ,
\end{gather*}
Therefore $(Y_n)_{n \geq 1}$ satisfies \eqref{cond2'} with $\nu=2$. Applying Theorem \ref{shaoresult}, we conclude that enlarging our probability space if necessary, there exists a sequence a sequence of independent centered Gaussian variables $(\widetilde Z_n)_{n \geq 1}$ such that
$\E(\widetilde Z_n^2)=\E(Y_n^2)$ and $Y_1+ \cdots Y_n-(\widetilde Z_1+\cdots +\widetilde Z_n)=o(n^{1/p} \sqrt{b(n) \log n})$ \as $\ $ We consider now the sequence of iid centered Gaussian variables $( Z_n)_{n \geq 1}$ with variance $\E(X_1^2)$ as defined in the proof of Corollary \ref{corFLIL}.  Notice then that
$$
\E (Z_k-\tilde Z_k)^2 = \big (\Vert X_k \Vert_2  - \Vert Y_k \Vert_2 \big )^2 \leq \Vert X_k  -  Y_k \Vert^2_2 \leq \E(X_k^2{\bf 1}_{\{|X_k|> k^{1/p}\}}) \, ,
$$
where for the last inequality, we have used the fact that  $\E (X_n | \G_{n+1}) = 0$ a.s. Hence by stationarity
\begin{equation*}
\sum_{n\ge 1} \E (Z_n-\tilde Z_n)^2/n^{2/p} \le
\E(X_1^2 \sum_{1\le n\le |X_1|^p} 1/n^{2/p})\le C \E(|X_1|^p)<\infty\, .
\end{equation*}
Therefore by the Kolmogorov theorem (or Lemma \ref{reverse}),
$\sum_{n\ge 1} (Z_n-\tilde Z_n)/n^{2/p}$ converges \as and by the Kronecker
lemma $Z_1+\ldots +Z_n-( \widetilde Z_1+\ldots +\widetilde Z_n)=
 o(n^{1/p}\sqrt{\log \log n})$
\as $ \ $
This achieves the proof of Corollary \ref{corrateASIP}. \hfill $\square$



\medskip

\subsection{Proof of Corollary \ref{corrateASIPp=4}.} 

The proof relies more deeply on the construction of Scott and Huggins \cite{HS}. We want to use Theorem \ref{shaoresult} without condition \eqref{cond2'}. Now the proof of Theorem \ref{shaoresult} relies
on Proposition \ref{huggins} and \eqref{cond2'} is used to ensure that condition \eqref{cond2} holds for an
auxiliary process. Instead of \eqref{cond2'} we will make use of Remark \ref{remlil}.
We define a reverse martingale $(R_n)_{n\ge 1}$,  by $R_n=\sum_{k\ge n} X_k/k$. Notice that $R_n$ is well defined in $\LL^2$ by Lemma \ref{reverse}.

For every $n\le -1$ define $\widetilde R_n:=R_{-n}$,  $\widetilde X_n:=X_{-n}$ and $\widetilde \G_n:=\G_{-n}$. Then $(\widetilde R_n,\widetilde \G_n)_{n\le -1}$ is a
martingale.

Enlarging our probability space if necessary,
we may consider a countable set of standard Brownian motions $(B^{(n)}_t)_{t\ge 0}$, $n\le -1$ that are independent of each others and of $(\widetilde X_{n})_{n\le -1}$. Notice that the process
$(\widetilde X_{n},(B^{(n)}_t)_{t\ge 0})_{n\le -1}$ with values in $\R\times \R^{\R^+}$, is stationary.

\smallskip

We now define a filtration $(\widetilde \h_t)_{t\le -1}$ as follows. For $n\le -1$ an integer, write
$\widetilde \h_n= \widetilde \G_n \vee \sigma\{ B^{(j)}_t, \, 0\le t <\infty ,\, - \infty < j\le n\}$. For every $t\le -1$, not an integer, write
$\widetilde \h_t=\widetilde \h_{[t]}\vee \{R_{[t]+1}+B^{([t]+1)}_{\phi(s)}, \, 0<s\le t-[t]\}$, where $[t]$ stands for the largest negative integer, not exceeding $t$, and $\phi$ is defined on $]0,1]$ by $\phi(s):=1/s-1$. Then, we define a continuous martingale with respect to $(\widetilde \h_t)_{t \leq -1}$
interpolating $(\widetilde R_n)$, by $\widetilde R_t=\E(\widetilde R_{[t]+1}|\widetilde \h_t)$, for every $t\le -1$.  Notice that
\begin{equation}\label{key}
\widetilde R_t=\widetilde R_{[t]} +\frac{\E (X_{-[t]-1}|\widetilde \h_t)}{-[t]-1} \, .
\end{equation}

Using Theorem A of \cite{HS} as done page 451 of \cite{HS}, there exists a continuous non-decreasing process
$(\widetilde \tau_t )_{t\le -1}$ and a Brownian motion $(B^*_t)_{t\ge 0}$ such that  $\widetilde R_t=B^*_{\widetilde \tau_t}$ a.s. and
$(\widetilde R_t^2-\widetilde \tau_t)_{t \leq -1}$ is a martingale with respect to $(\widetilde \h_t)_{t \leq -1}$. In particular, $(\widetilde \tau_t)_{t \leq -1}$ must be the quadratic variation of
$(\widetilde R_t)$ on $]-\infty, t]$.

For every $n\ge 1$, define $\tau_n:=\widetilde \tau_{-n}$, $\h_n:=\widetilde \h_{-n}$. These are exactly the quantities involved
in the proof of Proposition \ref{huggins}. Then $t_n=\tau_n-\tau_{n+1}$ is nothing else but the quadratic variation
of $(\widetilde R_t)$ on $[-n-1, -n]$.  But it follows from \eqref{key} that
$(n^2t_n)_{n \geq 1}$ is a stationary and ergodic process.

By Remark \ref{remlil} we need to prove
\begin{equation}\label{LL}
\sum_{k\ge n} (t_k-\E(t_k|\h_{k+1}))
=O \big (n^{-3/2}\sqrt {\log \log n} \big ) \quad \as
\end{equation}

\medskip

Since $(n^2t_n)_{n \geq 1}$ is a stationary and ergodic sequence in $\LL^2$, it follows from Corollary \ref{corFLIL} that
$\sum_{k=1}^n k^2(t_k-\E(t_k|\h_{k+1}))= O(\sqrt{n\log \log n})$ \as, which proves \eqref{LL} by an Abel summation.

\hfill $\square$

\end{document}